# Mordukhovich derivatives of the normalized duality mapping in Banach spaces

Jinlu Li

Department of Mathematics
Shawnee State University
Portsmouth, Ohio 45662 USA
jli@shawnee.edu

**Abstract**

In this paper, we investigate some properties of the Mordukhovich derivatives of the normalized duality mapping in Banach spaces. For the underlying spaces, we consider three cases: uniformly convex and uniformly smooth Banach space $l_p$; general Banach spaces $L_1$ and $C[0,1]$.

**Keywords**: Normalized duality mapping; Fréchet differentiability; Mordukhovich derivative; coderivative

## 1. Introduction

Let $(X, \|\cdot\|_X)$ and $(Y, \|\cdot\|_Y)$ be real Banach spaces with topological dual spaces $X^*$ and $Y^*$, respectively. Let $\langle \cdot, \cdot \rangle_X$ denote the real canonical pairing between $X^*$ and $X$ and $\langle \cdot, \cdot \rangle_Y$ the real canonical pairing between $Y^*$ and $Y$. Let $\Delta$ be a nonempty subset of $X$ and let $F: \Delta \rightrightarrows Y$ be a set valued mapping. The graph of $F$ is defined by the following subset in $\Delta \times Y$

$$\text{gph}F = \{(x,y) \in \Delta \times Y : y \in F(x)\}.$$

For $(x,y) \in \text{gph}F$, that is, for $x \in \Delta$ and $y \in F(x)$, the Mordukhovich derivative (which is also called Mordukhovich coderivative, or coderivative) of $F$ at point $(x, y)$ is a set valued mapping $\widehat{D}^*F(x,y): Y^* \rightrightarrows X^*$. For any $y^* \in Y^*$, it is defined by (see Definitions 1.13 and 1.32 in Chapter 1 in [20])

$$\widehat{D}^*F(x,y)(y^*) = \left\{ z^* \in X^* : \limsup_{\substack{(u,v) \to (x,y) \\ u \in \Delta \text{ and } v \in F(u)}} \frac{\langle z^*, u-x \rangle_X - \langle y^*, v-y \rangle_Y}{\|u-x\|_X + \|v-y\|_Y} \leq 0 \right\}$$

$$= \left\{ z^* \in X^* : \limsup_{\substack{(u,v) \to (x,y) \\ (u,v) \in \text{gph}F}} \frac{\langle z^*, u-x \rangle_X - \langle y^*, v-y \rangle_Y}{\|u-x\|_X + \|v-y\|_Y} \leq 0 \right\}. \qquad (1.1)$$

If $(x, y) \notin \text{gph}F$, then, we define

$$\widehat{D}^*F(x,y)(y^*) = \emptyset, \text{ for any } y^* \in Y^*.$$

By the above definition (1.1), $\widehat{D}^*F(x,y): Y^* \rightrightarrows X^*$ is a set valued mapping, which is called the Mordukhovich differential operator (or the Mordukhovich codifferential operator) of $F$ at $(x, y)$.

In particular, let $G: \Delta \to Y$ be a single valued continuous mapping. By (1.1), the Mordukhovich derivative of $G$ at point $(x, G(x))$ is a set valued mapping $\widehat{D}^*G(x, G(x)): Y^* \rightrightarrows X^*$. For any $y^* \in Y^*$, it is defined by

$$\widehat{D}^*G(x, G(x))(y^*) := \widehat{D}^*G(x)(y^*)$$

$$= \left\{ z^* \in X^*: \limsup_{\substack{(u,G(u)) \to (x,G(x)) \\ u \in \Delta}} \frac{\langle z^*, u-x \rangle_X - \langle y^*, G(u)-G(x) \rangle_Y}{\|u-x\|_X + \|G(u)-G(x)\|_Y} \leq 0 \right\}$$

$$= \left\{ z^* \in X^*: \limsup_{\substack{u \to x \\ u \in \Delta}} \frac{\langle z^*, u-x \rangle_X - \langle y^*, G(u)-G(x) \rangle_Y}{\|u-x\|_X + \|G(u)-G(x)\|_Y} \leq 0 \right\}. \qquad (1.2)$$

The Mordukhovich derivatives have been widely applied to several branches of mathematics such as operator theory, optimization theory, approximation theory, control theory, equilibrium theory, and so forth (see [20−22]).

For the single valued mapping $G$, the smoothness of $G$ is traditionally defined by some types of differentiability of $G$, such as Gâteaux directional differentiability, Fréchet differentiability and strict Fréchet differentiability. For example, see [9,12, 13, 23], in which the underlying spaces are Hilbert spaces; and see [4, 10, 26, 27], in which the underlying spaces are Banach spaces and normed linear spaces.

In contrast with the traditional differentiability of single valued mapping $G$, in [20], it is proved that if $G$ is Fréchet differentiable at a point $x$, then, the Mordukhovich derivative of $G$ at $x$ can be calculated in terms of the Fréchet derivative. Hence, the Mordukhovich derivatives can be considered as the generalization of Fréchet derivatives.

Since the metric projection operator is extremely important in operator theory, the present author studied the Mordukhovich derivatives of the metric projection operator in Hilbert spaces in [15,16]; in uniformly convex and uniformly smooth Banach spaces in [17]; and in general, in [18, 19].

Except the metric projection operator, it is well known that the normalized duality mapping $J$ is one of the most important mappings in the analysis in Banach spaces. It has many useful properties, which have been widely applied to fixed point theory, optimization theory, variational analysis, approximation theory, and so forth (see [1, 5, 20, 21, 22]). In this paper, we study the Mordukhovich differentiability of the normalized duality mapping in Banach spaces.

If the underlying Banach spaces are uniformly convex and uniformly smooth, then the normalized duality mapping $J$ is a single valued continuous and onto mapping. In section 3, by definition (1.2), we consider the Mordukhovich differentiability of $J$ in uniformly convex and uniformly smooth Banach space $l_p$, for $1 < p < \infty$.

In general, the normalized duality mapping $J$ in Banach spaces is a set valued mapping. In section 4, we use definition (1.1) to study the Mordukhovich differentiability of $J$ in general Banach space $L_1$. In section 5, we study the Mordukhovich differentiability of $J$ in $C[0, 1]$.

If we compare the results in this paper and the results in [20−22], we can find the significant differences between the Mordukhovich derivatives of the normalized duality mapping and the Mordukhovich derivatives of the metric projection operator, with respect to the same underlying Banach spaces.

2. **Preliminaries**

Let $(X, \|\cdot\|)$ be a general (real) Banach space with topological dual space $(X^*, \|\cdot\|_*)$. The dual space of $(X^*, \|\cdot\|_*)$ is denoted by $(X^{**}, \|\cdot\|_{**})$. Let $\langle \cdot, \cdot \rangle$ denote the real canonical pairing between $X^*$ and $X$; and let $\langle \cdot, \cdot \rangle_*$ denote the real canonical pairing between $X^{**}$ and $X^*$. Let $\theta$, $\theta^*$ and $\theta^{**}$ denote the origins in $X$, $X^*$ and $X^{**}$, respectively. The identity mappings on $X$, $X^*$ and $X^{**}$ are respectively denoted by $I_X$, $I_{X^*}$ and $I_{X^{**}}$.

The normalized duality mapping $J: X \rightrightarrows X^*$ is defined by

$$J(x) \text{ (or } Jx) = \{jx \in X^*: \langle jx, x \rangle = \|jx\|_* \|x\| = \|x\|^2 = \|jx\|_*^2\}, \text{ for any } x \in X.$$

The normalized duality mapping on $X^*$ is similarly denoted by $J^*: X^* \rightrightarrows X^{**}$, which is defined by, for any $x^* \in X^*$,

$$J^*(x^*) = \{j^*(x^*) \in X^{**}: \langle j^*(x^*), x^* \rangle_* = \|j^*(x^*)\|_{**} \|x^*\|_* = \|j^*(x^*)\|_{**}^2 = \|x^*\|_*^2\}.$$

When the considered Banach space is not reflexive, for $x^* \in X^*$, the value $J^*(x^*)$ contains a subset in $X$.

For $(x, x^*) \in \mathrm{gph} J$, that is, for $x \in X$ and $x^* \in J(x)$, by definition (1.1), the Mordukhovich derivative of $J$ at point $(x, x^*)$ is a set valued mapping $\widehat{D}^* J(x, x^*): X^{**} \rightrightarrows X^*$. For any $y^{**} \in X^{**}$, it is defined by

$$\widehat{D}^* J(x, x^*)(y^{**}) = \left\{ z^* \in X^*: \limsup_{\substack{(u,u^*) \to (x,x^*) \\ u \in X \text{ and } u^* \in J(u)}} \frac{\langle z^*, u-x \rangle - \langle y^{**}, u^*-x^* \rangle_*}{\|u-x\| + \|u^*-x^*\|_*} \leq 0 \right\}$$

$$= \left\{ z^* \in X^*: \limsup_{\substack{(u,u^*) \to (x,x^*) \\ (u,u^*) \in \mathrm{gph} J}} \frac{\langle z^*, u-x \rangle - \langle y^{**}, u^*-x^* \rangle_*}{\|u-x\| + \|u^*-x^*\|_*} \leq 0 \right\}. \quad (2.1)$$

In particular, if $X$ is a uniformly convex and uniformly smooth Banach space, then, $X$ is reflexive with $X^{**} = X$. In this case, it is well known that the normalized duality mapping $J: X \to X^*$ is a single valued continuous mapping. By (1.2), the Mordukhovich derivative of $J$ at point $(x, J(x))$ is a set valued mapping $\widehat{D}^*J(x,J(x)) : X \rightrightarrows X^*$ that is defined by, for any $y \in X$,

$$\widehat{D}^*J(x,J(x))(y) := \widehat{D}^*J(x)(y)$$

$$= \left\{ z^* \in X^*: \limsup_{\substack{(u,J(u))\to(x,J(x)) \\ u\in X}} \frac{\langle z^*,u-x\rangle - \langle y,J(u)-J(x)\rangle_*}{\|u-x\| + \|J(u)-J(x)\|_*} \leq 0 \right\}$$

$$= \left\{ z^* \in X^*: \limsup_{u\to x} \frac{\langle z^*,u-x\rangle - \langle J(u)-J(x), y\rangle}{\|u-x\| + \|J(u)-J(x)\|_*} \leq 0 \right\}. \tag{2.2}$$

Since the theme of this paper is about the normalized duality mapping in Banach spaces, we list some properties of the normalized duality mapping in the Appendix for easy reference. For more details, one may see Sections 4.2–4.3 and Problem set 4.2 in [28] and [1−3, 7, 11,14, 24, 25].

### 3. The normalized duality mapping in uniformly convex and uniformly smooth Banach space $l_p$, for $1 < p < \infty$

In this section, we focus on the real uniformly convex and uniformly smooth Banach spaces $(l_p, \|\cdot\|_p)$ and $(l_q, \|\cdot\|_q)$, in which $p$ and $q$ satisfy $1 < p, q < \infty$ and $\frac{1}{p} + \frac{1}{q} = 1$. The spaces $(l_p, \|\cdot\|_p)$ and $(l_q, \|\cdot\|_q)$ are the dual spaces of each other. That is, $l_p^* = l_q$ and $l_p^{**} = l_q^* = l_p$. Both $l_p$ and $l_q$ have origin $\theta = (0, 0, \ldots)$. Let $\langle \cdot, \cdot \rangle$ denote the real canonical pairing between $l_q$ and $l_p$; and let $\langle \cdot, \cdot \rangle_*$ denote the real canonical pairing between $l_p$ and $l_q$. We use $x, y, z, \ldots$ for the elements in $l_p$ and $u, v, w, \ldots$ for the elements in $l_q$.

The positive cone of $l_p$ is denoted by $l_p^+$ that is defined by

$$l_p^+ = \{x = (x_1, x_2, \ldots) \in l_p : x_n \geq 0, \text{ for all } n\}.$$

$l_p^+$ is a pointed closed and convex cone in $l_p$. We define a subset $l_p^{++}$ in $l_p$ as follows

$$l_p^{++} = \{x = (x_1, x_2, \ldots) \in l_p : x_n > 0, \text{ for all } n\}.$$

$l_p^{++}$ is a convex subset in $l_p$. However, $\{\theta\} \cup l_p^{++}$ is a pointed convex cone in $l_p$ with vertex $\theta$, which is neither closed, nor open.

The normalized duality mapping on $l_p$ is denoted by $J: l_p \to l_q$ and the normalized duality mapping on $l_p^* = l_q$ is denoted by $J^*: l_q \to l_p$. Recall the representations of normalized duality mapping $J: l_p \to l_q$ and the normalized duality mapping $J^*: l_q \to l_p$. For any point $x = (x_1, x_2, \ldots) \in l_p$ with $x \neq \theta$, we have

$$J(x) = \left( \frac{|x_1|^{p-1}\text{sign}(x_1)}{\|x\|_p^{p-2}}, \frac{|x_2|^{p-1}\text{sign}(x_2)}{\|x\|_p^{p-2}}, \ldots \right) = \left( \frac{|x_1|^{p-2}x_1}{\|x\|_p^{p-2}}, \frac{|x_2|^{p-2}x_2}{\|x\|_p^{p-2}}, \ldots \right). \tag{3.1}$$

Similar to (3.1), for any $u = (u_1, u_2, \ldots) \in l_q$ with $u \neq \theta$, we have

$$J^*(u) = \left(\frac{|u_1|^{q-1}\text{sign}(u_1)}{\|u\|_q^{q-2}}, \frac{|u_2|^{q-1}\text{sign}(u_2)}{\|u\|_q^{q-2}}, \ldots\right) = \left(\frac{|u_1|^{q-2}u_1}{\|u\|_q^{q-2}}, \frac{|u_2|^{q-2}u_2}{\|u\|_q^{q-2}}, \ldots\right). \quad (3.2)$$

By (3.1), $J$ is a mapping from $l_p^+$ to $l_q^+$; and from $l_p^{++}$ to $l_q^{++}$. That is,

$$x \in l_p^+ \implies J(x) \in l_q^+ \quad \text{and} \quad x \in l_p^{++} \implies J(x) \in l_q^{++}.$$

It follows that, for any point $x = (x_1, x_2, \ldots) \in l_p^+$ with $x \neq \theta$, we have

$$J(x) = \left(\frac{x_1^{p-1}}{\|x\|_p^{p-2}}, \frac{x_2^{p-1}}{\|x\|_p^{p-2}}, \ldots\right).$$

Similar to (3.1), for any $u = (u_1, u_2, \ldots) \in l_q^+$ with $u \neq \theta$, we have

$$J^*(u) = \left(\frac{u_1^{q-1}}{\|u\|_q^{q-2}}, \frac{u_2^{q-1}}{\|u\|_q^{q-2}}, \ldots\right).$$

One has $J^* \circ J = I_{l_p}$ and $J \circ J^* = I_{l_q}$. Since both mappings $J: l_p \to l_q$ and $J^*: l_q \to l_p$ are single valued continuous mappings, for any $x \in l_p$ with $J(x) \in l_q$, by (2.2), the Mordukhovich derivative of $J$ at $(x, J(x))$ is a set valued mapping $\widehat{D}^*J(x, J(x)): l_q^* \rightrightarrows l_p^*$. It can be simply rewritten as, $\widehat{D}^*J(x) := \widehat{D}^*J(x, J(x)): l_p \rightrightarrows l_q$. By (2.2), for any $y \in l_p = l_q^*$, it is defined by

$$\widehat{D}^*J(x, J(x))(y) := \widehat{D}^*J(x)(y)$$

$$= \left\{w \in l_q : \limsup_{z \to x} \frac{\langle w, z-x \rangle - \langle y, J(z) - J(x) \rangle_*}{\|z-x\|_p + \|J(z) - J(x)\|_q} \leq 0\right\}$$

$$= \left\{w \in l_q : \limsup_{z \to x} \frac{\langle w, z-x \rangle - \langle J(z) - J(x), y \rangle}{\|z-x\|_p + \|J(z) - J(x)\|_q} \leq 0\right\}. \quad (3.3)$$

**Theorem 3.1.** *Let* $x = (x_1, x_2, \ldots) \in l_p$. *Then*

$$\widehat{D}^*J(x)(\theta) = \{\theta\}.$$

*Proof.* It is clear to see that

$$\theta \in \widehat{D}^*J(x)(\theta). \quad (3.4)$$

Next, we prove that, for any $w \in l_q$, if $w \neq \theta$, then

$$w \notin \widehat{D}^*J(x)(\theta). \quad (3.5)$$

The proof of (3.5) is divided into two cases with respect to $x$.

Case 1. $x \neq \theta$. Let $w = (w_1, w_2, \ldots) \in l_q$ with $w \neq \theta$. There is a positive integer $m$ such that $w_m \neq 0$. We may assume $w_m > 0$. Let $\lambda_m \in l_p \cap l_q$, in which, the $m$th coordinate is 1 and all other coordinates are 0. In this case, in the limit in (3.3), we take a special direction $z_t = t\lambda_m + x$, for $t > 0$ with $t \downarrow 0$. More precisely speaking, for all $n$, the $n$th coordinate of $z_t$ has the following representation.

$$(z_t)_n = \begin{cases} t + x_m, & \text{for } n = m, \\ x_n, & \text{for } n \neq m. \end{cases}$$

This implies that $\|z_t - x\|_p = t \to 0$, as $t \downarrow 0$. That is,

$$z_t \to x \text{ in } l_p, \text{ as } t \downarrow 0. \tag{3.6}$$

For $t > 0$, by $z_t = t\lambda_m + x$ and $x = z_t - t\lambda_m$, we have

$$\|x\|_p - t < \|z_t\|_p < \|x\|_p + t, \text{ for all } t > 0.$$

This implies

$$\big|\|z_t\|_p - \|x\|_p\big| < t, \text{ for all } t > 0. \tag{3.7}$$

In particular,

$$\frac{1}{2}\|x\|_p < \|z_t\|_p < \frac{3}{2}\|x\|_p, \text{ for all } 0 < t < \frac{1}{2}\|x\|_p. \tag{3.8}$$

By (3.6) and by the continuity of the normalized duality mapping in uniformly convex and uniformly smooth Banach spaces, we have

$$J(z_t) \to J(x) \text{ in } l_q, \text{ as } t \downarrow 0.$$

However, we need more estimation for $\|J(z_t) - J(x)\|_q$. By (3.1), we have

$$J(x) = \left( \frac{|x_1|^{p-1}\text{sign}(x_1)}{\|x\|_p^{p-2}}, \frac{|x_2|^{p-1}\text{sign}(x_2)}{\|x\|_p^{p-2}}, \ldots \frac{|x_m|^{p-1}\text{sign}(x_m)}{\|x\|_p^{p-2}}, \ldots \right)$$

and

$$J(z_t) = \left( \frac{|x_1|^{p-1}\text{sign}(x_1)}{\|z_t\|_p^{p-2}}, \frac{|x_2|^{p-1}\text{sign}(x_2)}{\|z_t\|_p^{p-2}}, \ldots \frac{|t+x_m|^{p-1}\text{sign}(t+x_m)}{\|z_t\|_p^{p-2}}, \ldots \right).$$

It follows that

$$J(z_t) - J(x) = \left( \frac{1}{\|z_t\|_p^{p-2}} - \frac{1}{\|x\|_p^{p-2}} \right)(\|x\|_p^{p-2})J(x) + \left( \frac{|t+x_m|^{p-1}\text{sign}(t+x_m)}{\|z_t\|_p^{p-2}} - \frac{|x_m|^{p-1}\text{sign}(x_m)}{\|z_t\|_p^{p-2}} \right)\lambda_m.$$

This implies

$$\|J(z_t) - J(x)\|_q$$

$$\leq \left\| \left( \frac{1}{\|z_t\|_p^{p-2}} - \frac{1}{\|x\|_p^{p-2}} \right)(\|x\|_p^{p-2})J(x) \right\|_q + \left\| \left( \frac{|t+x_m|^{p-1}\text{sign}(t+x_m)}{\|z_t\|_p^{p-2}} - \frac{|x_m|^{p-1}\text{sign}(x_m)}{\|z_t\|_p^{p-2}} \right)\lambda_m \right\|_q$$

$$= \left| \frac{1}{\|z_t\|_p^{p-2}} - \frac{1}{\|x\|_p^{p-2}} \right| \|x\|_p^{p-2}\|J(x)\|_q + \left| \frac{|t+x_m|^{p-1}\text{sign}(t+x_m)}{\|z_t\|_p^{p-2}} - \frac{|x_m|^{p-1}\text{sign}(x_m)}{\|z_t\|_p^{p-2}} \right|$$

$$= \left| \frac{1}{\|z_t\|_p^{p-2}} - \frac{1}{\|x\|_p^{p-2}} \right| \|x\|_p^{p-1} + \frac{1}{\|z_t\|_p^{p-2}} \left||t+x_m|^{p-1}\text{sign}(t+x_m) - |x_m|^{p-1}\text{sign}(x_m)\right|. \quad (3.9)$$

For the estimation of $\frac{1}{\|z_t\|_p^{p-2}}\left||t+x_m|^{p-1}\text{sign}(t+x_m) - |x_m|^{p-1}\text{sign}(x_m)\right|$, we consider 3 subcases below with respect to $x_m$.

Subcase 1.1. $x_m = 0$. In this subcase, by $p - 1 > 0$, we have

$$\left||t+x_m|^{p-1}\text{sign}(t+x_m) - |x_m|^{p-1}\text{sign}(x_m)\right| = t^{p-1}, \text{ for all } t > 0.$$

Since $x \neq \theta$ and $x_m = 0$, we have $\|z_t\|_p > t$. This implies

$$\frac{1}{\|z_t\|_p^{p-2}}\left||t+x_m|^{p-1}\text{sign}(t+x_m) - |x_m|^{p-1}\text{sign}(x_m)\right|$$

$$= \frac{t^{p-1}}{\|z_t\|_p^{p-2}} < \frac{t^{p-1}}{t^{p-2}} = t. \quad (3.10)$$

Subcase 1.2. $x_m > 0$. In this subcase, for all $1 > t > 0$, by $p - 1 > 0$, we have

$$\left||t+x_m|^{p-1}\text{sign}(t+x_m) - |x_m|^{p-1}\text{sign}(x_m)\right|$$

$$= (t+x_m)^{p-1} - x_m^{p-1}$$

$$= (p-1)\,\bar{x}^{p-2}(t + x_m - x_m)$$

$$= (p-1)t\,\bar{x}^{p-2}, \text{ for some } \bar{x} \in (x_m, x_m + t). \quad (3.11)$$

In this subcase, for all $t > 0$, we have $\|z_t\|_p > \|x\|_p$. By (3.11), for all $1 > t > 0$, this implies,

$$\frac{1}{\|z_t\|_p^{p-2}}\left||t+x_m|^{p-1}\text{sign}(t+x_m) - |x_m|^{p-1}\text{sign}(x_m)\right|$$

$$= \frac{(p-1)t\,\bar{x}^{p-2}}{\|z_t\|_p^{p-2}}$$

$$< \frac{(p-1)\,\bar{x}^{p-2}}{\|x\|_p^{p-2}}t, \text{ for some } \bar{x} \in (x_m, x_m + 1).$$

$$\leq \frac{(p-1)(x_m+\beta)^{p-2}}{\|x\|_p^{p-2}}t. \quad (3.12)$$

Where, we take $\beta = 1$, for $p - 2 \geq 0$ and $\beta = 0$, for $p - 2 < 0$.

Subcase 1.3. $x_m < 0$. In this case, if $t < -x_m$, then $t + x_m < 0$. Then, for all $0 < t < -\frac{1}{2}x_m$, by (3.7), we have $\|z_t\|_p > \frac{1}{2}\|x\|_p$. This implies

$$\frac{1}{\|z_t\|_p^{p-2}}\left||t+x_m|^{p-1}\text{sign}(t+x_m) - |x_m|^{p-1}\text{sign}(x_m)\right|$$

$$= \frac{1}{\|z_t\|_p^{p-2}} \left| -(-(t+x_m))^{p-1} - (-(-x_m)^{p-1}) \right|$$

$$= \frac{1}{\|z_t\|_p^{p-2}} \left| (-(t+x_m))^{p-1} - (-x_m)^{p-1} \right|$$

$$= \frac{1}{\|z_t\|_p^{p-2}} \left( (-x_m)^{p-1} - (-(t+x_m))^{p-1} \right)$$

$$= \frac{1}{\|z_t\|_p^{p-2}} (p-1)\, \tilde{x}^{p-2} (-x_m - (-(t+x_m)))$$

$$= \frac{(p-1)t\, \tilde{x}^{p-2}}{\|z_t\|_p^{p-2}}$$

$$< \frac{2(p-1)\, \tilde{x}^{p-2}}{\|x\|_p^{p-2}} t, \text{ for some } -(t+x_m) < \tilde{x} < -x_m$$

$$\leq \frac{2(p-1)\, (-x_m+\alpha)^{p-2}}{\|x\|_p^{p-2}} t. \tag{3.13}$$

Where, we take $\alpha = -\frac{1}{2} x_m$, for $p - 2 \geq 0$ and $\alpha = 0$, for $p - 2 < 0$. From (3.10), (3.12) and (3.13), let

$$A = \max\left\{1, \frac{(p-1)\,(x_m+\beta)^{p-2}}{\|x\|_p^{p-2}}, \frac{(p-1)\,(-x_m+\alpha)^{p-2}}{\|x\|_p^{p-2}}\right\}.$$

Then, we have

$$\frac{1}{\|z_t\|_p^{p-2}} ||t+x_m|^{p-1}\operatorname{sign}(t+x_m) - |x_m|^{p-1}\operatorname{sign}(x_m)| < At, \text{ as } 0 < t < \frac{1}{2}|x_m|. \tag{3.14}$$

Next, we estimate $\left| \frac{1}{\|z_t\|_p^{p-2}} - \frac{1}{\|x\|_p^{p-2}} \right|$, which is the first term in (3.9). By (3.7) and (3.8), we have

$$\left| \frac{1}{\|z_t\|_p^{p-2}} - \frac{1}{\|x\|_p^{p-2}} \right|$$

$$= \frac{\left| \|z_t\|_p^2 \|x\|_p^p - \|x\|_p^2 \|z_t\|_p^p \right|}{\|z_t\|_p^p \|x\|_p^p}$$

$$= \frac{\left| \|z_t\|_p^2 \|x\|_p^p - \|x\|_p^2 \left( \|x\|_p^p + |t+x_m|^p - |x_m|^p \right) \right|}{\|z_t\|_p^p \|x\|_p^p}$$

$$\leq \frac{\left| \|z_t\|_p^2 \|x\|_p^p - \|x\|_p^2 \|x\|_p^p \right|}{\|z_t\|_p^p \|x\|_p^p} + \frac{\|x\|_p^2 ||t+x_m|^p - |x_m|^p|}{\|z_t\|_p^p \|x\|_p^p}$$

$$= \frac{\left| \|z_t\|_p^2 - \|x\|_p^2 \right|}{\|z_t\|_p^p} + \frac{||t+x_m|^p - |x_m|^p|}{\|z_t\|_p^p \|x\|_p^{p-2}}$$

$$= \frac{(\|z_t\|_p + \|x\|_p)\|\|z_t\|_p - \|x\|_p\|}{\|z_t\|_p^p} + \frac{\|\|t + x_m\|^p - \|x_m\|^p\|}{\|z_t\|_p^p \|x\|_p^{p-2}}$$

$$< \frac{\frac{5}{2}\|x\|_p \, t}{\frac{1}{2}\|x\|_p^p} + \frac{\|\|t+x_m\|^p - \|x_m\|^p\|}{\frac{1}{2^p}\|x\|_p^{2p-2}}$$

$$= \frac{5\,t}{\|x\|_p^{p-1}} + \frac{2^p \|\|t+x_m\|^p - \|x_m\|^p\|}{\|x\|_p^{2p-2}}. \tag{3.15}$$

Similar to the proof of (3.14), we can prove that there is a positive number $B$ such that

$$\|\|t + x_m\|^p - \|x_m\|^p\| < Bt, \text{ as } 0 < t < \frac{1}{2}|x_m|. \tag{3.16}$$

By (3.15) and (3.16), we obtain

$$\left| \frac{1}{\|z_t\|_p^{p-2}} - \frac{1}{\|x\|_p^{p-2}} \right| < \frac{5\,t}{\|x\|_p^{p-1}} + \frac{2^p Bt}{\|x\|_p^{2p-2}}, \text{ as } 0 < t < \frac{1}{2}|x_m|. \tag{3.17}$$

Substituting the results of (3.17) and (3.14) into (3.9), we get

$$\|J(z_t) - J(x)\|_q$$

$$< \left( \frac{5\,t}{\|x\|_p^{p-1}} + \frac{2^p Bt}{\|x\|_p^{2p-2}} \right) \|x\|_p^{p-1} + At$$

$$= 5\,t + \frac{2^p Bt}{\|x\|_p^{p-1}} + At$$

$$:= Ct, \text{ as } 0 < t < \frac{1}{2}|x_m|. \tag{3.18}$$

Where $C = 5\,t + \frac{2^p Bt}{\|x\|_p^{p-1}} + At$ satisfying $C > 0$, for $0 < t < \frac{1}{2}|x_m|$. This implies

$$J(z_t) \to J(x) \text{ in } l_q, \text{ as } t \downarrow 0. \tag{3.19}$$

By (3.6) and (3.19), we have

$$(z_t, J(z_t)) \to (x, J(x)), \text{ as } t \downarrow 0. \tag{3.20}$$

In order to show $w \notin \widehat{D}^* J(x)(\theta)$, we calculate the limit in (3.1). By the assumption that $w_m > 0$ and by (3.20) and (3.18), we have

$$\limsup_{z \to x} \frac{\langle w, z-x \rangle - \langle J(z) - J(x), \theta \rangle}{\|z-x\|_p + \|J(z) - J(x)\|_q}$$

$$= \limsup_{z \to x} \frac{\langle w, z-x \rangle}{\|z-x\|_p + \|J(z) - J(x)\|_q}$$

$$\geq \limsup_{z_t \to x} \frac{\langle w, z_t - x \rangle}{\|z_t - x\|_p + \|J(z_t) - J(x)\|_q}$$

$$= \limsup_{t \downarrow 0} \frac{\langle w, t\lambda_m \rangle}{\|t\lambda_m\|_p + \|J(z_t) - J(x)\|_q}$$

$$= \limsup_{t \downarrow 0} \frac{tw_m}{t + \|J(z_t) - J(x)\|_q}$$

$$> \limsup_{t \downarrow 0} \frac{tw_m}{t + Ct}$$

$$= \frac{w_m}{1 + C} > 0.$$

This implies that, for $x \in l_p$ with $x \neq \theta$, we have

$$w \notin \widehat{D}^*J(x)(\theta), \text{ for any } w \in l_p \text{ with } w \neq \theta. \tag{3.21}$$

By (3.4) and (3.21), we obtain

$$\widehat{D}^*J(x)(\theta) = \{\theta\}, \text{ for } x \in l_p \text{ with } x \neq \theta. \tag{3.22}$$

Case 2. $x = \theta$, For the case that $x = \theta$, as in the proof of (3.22), let $w \in l_p$ with $w_m > 0$, for some positive integer $m$. For $t > 0$, let $z_t = t\lambda_m + \theta = t\lambda_m$, we have

$$\big|\|z_t\|_p - \|\theta\|_p\big| = t, \text{ for all } t > 0.$$

and
$$J(z_t) = J(t\lambda_m) = t\lambda_m \to \theta, \text{ as } t \downarrow 0.$$

By $J\theta = \theta$, these imply

$$(z_t, J(z_t)) = (t\lambda_m, t\lambda_m) \to (\theta, \theta), \text{ as } t \downarrow 0.$$

We have

$$\limsup_{z \to \theta} \frac{\langle w, z - \theta \rangle - \langle J(z) - J(\theta), \theta \rangle}{\|z - \theta\|_p + \|J(z) - J(\theta)\|_q}$$

$$\geq \limsup_{z_t \to \theta} \frac{\langle w, z_t - \theta \rangle}{\|z_t - \theta\|_p + \|J(z_t) - J(\theta)\|_q}$$

$$= \limsup_{t \downarrow 0} \frac{\langle w, t\lambda_m \rangle}{\|t\lambda_m\|_p + \|t\lambda_m\|_q}$$

$$= \limsup_{t \downarrow 0} \frac{tw_m}{2t}$$

$$= \frac{w_m}{2} > 0.$$

This implies

$$w \notin \widehat{D}^*J(\theta)(\theta), \text{ for any } w \in l_p \text{ with } w \neq \theta.$$

By (3.4), (3.22), this proves this theorem. □

**Theorem 3.2.** *Let* $x, y \in l_p$. *If* $\langle J(x), y \rangle \neq 0$, *then*

$$\theta \notin \widehat{D}^*J(x)(y).$$

*Proof.* For $t > 0$, let

$$z_t = \begin{cases} (1-t)x, & \text{if } \langle J(x), y \rangle > 0, \\ (1+t)x, & \text{if } \langle J(x), y \rangle < 0. \end{cases}$$

By the properties of the normalized duality mapping, for $t > 0$, we have

$$J(z_t) = \begin{cases} (1-t)J(x), & \text{if } \langle J(x), y \rangle > 0, \\ (1+t)J(x), & \text{if } \langle J(x), y \rangle < 0. \end{cases}$$

This implies

$$(z_t, J(z_t)) \to (x, J(x)), \text{ as } t \downarrow 0.$$

In order to show $\theta \notin \widehat{D}^*J(x)(y)$, we calculate the limit in (3.1).

$$\limsup_{z \to x} \frac{\langle \theta, z-x \rangle - \langle y, J(z)-J(x) \rangle_*}{\|z-x\|_p + \|J(z)-J(x)\|_q}$$

$$= \limsup_{z \to x} \frac{\langle \theta, z-x \rangle - \langle J(z)-J(x), y \rangle}{\|z-x\|_p + \|J(z)-J(x)\|_q}$$

$$= \limsup_{z \to x} \frac{-\langle J(z)-J(x), y \rangle}{\|z-x\|_p + \|J(z)-J(x)\|_q}$$

$$\geq \limsup_{z_t \to x} \frac{-\langle J(z_t)-J(x), y \rangle}{\|tx\|_p + \|tJ(x)\|_q}$$

$$= \limsup_{t \downarrow 0} \frac{t|\langle J(x), y \rangle|}{2t\|x\|_p}$$

$$= \frac{|\langle J(x), y \rangle|}{2\|x\|_p} > 0.$$

This implies

$$\theta \notin \widehat{D}^*J(x)(y), \text{ if } \langle J(x), y \rangle \neq 0.$$ □

**Theorem 3.3.** *Let* $x \in l_p \backslash \{\theta\}$. *For any* $a > 0$ *with* $a \neq 1$, *we have*

$$aJ(x) \notin \widehat{D}^*J(x)(x).$$

*Proof.* For $1 > t > 0$, let

$$z_t = \begin{cases} (1+t)x, & \text{if } a > 1, \\ (1-t)x, & \text{if } 0 < a < 1. \end{cases}$$

By the properties of the normalized duality mapping, for $1 > t > 0$, we have

$$J(z_t) = \begin{cases} (1+t)J(x), & \text{if } a > 1, \\ (1-t)J(x), & \text{if } 0 < a < 1. \end{cases}$$

This implies

$$(z_t, J(z_t)) \to (x, J(x)), \text{ as } t \downarrow 0.$$

In order to show $aJ(x) \notin \widehat{D}^*J(x, J(x))(x)$, we calculate the limit in (3.1).

$$\limsup_{z \to x} \frac{\langle aJ(x), z-x \rangle - \langle x, J(z)-J(x) \rangle_*}{\|z-x\|_p + \|J(z)-J(x)\|_q}$$

$$= \limsup_{z \to x} \frac{\langle aJ(x), z-x \rangle - \langle J(z)-J(x), x \rangle}{\|z-x\|_p + \|J(z)-J(x)\|_q}$$

$$\geq \limsup_{z_t \to x} \frac{\langle aJ(x), z_t-x \rangle - \langle J(z_t)-J(x), x \rangle}{\|tx\|_p + \|tJ(x)\|_q}$$

$$= \begin{cases} \limsup_{z_t \to x} \frac{\langle aJ(x), tx \rangle - \langle tJ(x), x \rangle}{\|tx\|_p + \|tJ(x)\|_q}, & a > 1 \\ \limsup_{z_t \to x} \frac{\langle aJ(x), -tx \rangle - \langle -tJ(x), x \rangle}{\|tx\|_p + \|tJ(x)\|_q}, & 0 < a < 1 \end{cases}$$

$$= \begin{cases} \limsup_{z_t \to x} \frac{t(a-1)\|x\|_p^2}{\|tx\|_p + \|tJ(x)\|_q}, & a > 1 \\ \limsup_{z_t \to x} \frac{t(1-a)\|x\|_p^2}{\|tx\|_p + \|tJ(x)\|_q}, & 0 < a < 1 \end{cases}$$

$$= \limsup_{t \downarrow 0} \frac{t|a-1|\|x\|_p^2}{2t\|x\|_p}$$

$$= \frac{|a-1|\|x\|_p}{2} > 0.$$

This proves this Theorem. □

## 4. The normalized duality mapping in (general) Banach space $L_1(S)$

Let $(S, \mathcal{A}, \mu)$ be a positive and complete measure space. The real Banach space $(L_1(S), \|\cdot\|_1)$ has dual space $(L_\infty(S), \|\cdot\|_\infty)$. Both $L_1(S)$ and $L_\infty(S)$ are not reflexive. We define the positive cone, denoted by $L_1^+(S)$, in $L_1(S)$ as follows.

$$L_1^+(S) = \{f \in L_1(S) : f(s) \geq 0, \text{ for almost all } s \in S\}.$$

It is well-known that $L_1^+(S)$ is a pointed closed and convex cone in $L_1(S)$. Then, we define a

subset in $L_1(S)$ as follows

$$L_1^{++}(S) = \{f \in L_1(S): f(s) > 0, \text{ for almost all } s \in S\}.$$

$L_1^{++}(S)$ is a convex subset in $L_1(S)$. However, $\{\theta\} \cup L_1^{++}(S)$ is a pointed convex cone in $L_1(S)$ with vertex $\theta$, which is neither closed, nor open.

Let $ba(S, \mathcal{A})$ be the Banach space of all bounded finitely additive real valued functions on $\mathcal{A}$ with norm $\|\cdot\|$. For every $\gamma \in ba(S, \mathcal{A})$, $\|\gamma\|$ is defined by (see page 160 of section III 7 in [8] by Dunford and Schwartz)

$$\|\gamma\| = \sup\{|\gamma(A)|: A \in \mathcal{A}\}, \text{ for every } \gamma \in ba(S, \mathcal{A}).$$

By Theorem IV 8.1 in [8] by Dunford and Schwartz, the Banach space $L_\infty(S)$ has dual space $L_\infty^*(S) = ba(S, \mathcal{A})$. Let $\langle \cdot, \cdot \rangle$ denote the real canonical pairing between $L_1(S)$ and $L_\infty(S)$; and let $\langle \cdot, \cdot \rangle_*$ denote the real canonical pairing between $L_\infty(S)$ and $ba(S, \mathcal{A})$. It satisfies that, $\varphi \in L_\infty^*(S)$ if and only if, there is $\gamma \in ba(S, \mathcal{A})$ such that

$$\langle \varphi, h^* \rangle_* = \int_S h^*(s)\gamma(ds), \text{ for any } h^* \in L_\infty(S) = L_1^*(S). \tag{4.1}$$

Then, $\varphi$ and $\gamma$ satisfying (4.1) are identified; and therefore, we have

$$L_1^*(S) = L_\infty(S) \text{ and } L_1^{**}(S) = L_\infty^*(S) = ba(S, \mathcal{A}).$$

We use English letters, such as $f, g, \ldots$, to name the elements in $L_1(S)$, and $f^*, g^*, \ldots$, for elements in $L_\infty(S)$. We use Greek letters, such as $\gamma, \lambda, \ldots$ to name the elements in $ba(S, \mathcal{A})$. Let $\theta$, $\theta^*$ and $\theta^{**}$ denote the origins in $L_1(S)$, $L_\infty(S)$ and $ba(S, \mathcal{A})$, respectively. The following lemma shows that the normalized duality mapping $J: L_1(S) \rightrightarrows L_\infty(S)$ is indeed a set valued mapping.

**Lemma 4.1.** *Let $f \in L_1(S)$ with $f \neq \theta$. Let a be an arbitrarily given bounded measurable function defined on the set $\{s \in S: f(s) = 0\}$ satisfying*

$$-\|f\|_1 \leq a(s) \leq \|f\|_1, \text{ for all } s \in S \text{ with } f(s) = 0.$$

*Define jf $\in L_\infty(S)$ as follows*

$$(jf)(s) = \begin{cases} \|f\|_1, & \text{for } f(s) > 0, \\ a(s), & \text{for } f(s) = 0, \text{ for all } s \in S. \\ -\|f\|_1, & \text{for } f(s) < 0, \end{cases} \tag{4.2}$$

*Then, jf $\in$ Jf.*

*Proof.* Let $A = \{s \in S: f(s) > 0\}$ and $B = \{s \in S: f(s) < 0\}$. Since $f \neq \theta$, then $\mu(A \cup B) \neq 0$. By (4.2), this implies $\|jf\|_\infty = \|f\|_1 > 0$. Then,

$$\langle jf, f \rangle = \int_S f(s)(jf)(s)\mu(ds)$$

$$= \|f\|_1 \int_S |f(s)| \mu(ds)$$

$$= \|f\|_1^2 = \|jf\|_\infty^2.$$

This proves $jf \in Jf$. □

**Corollary 4.2.** *Let $f \in L_1(S)$ with $f \neq 0$. We have*

(i) *If $\mu\{s \in S: f(s) = 0\} > 0$, then $Jf$ is an infinite set;*
(ii) *If $\mu\{s \in S: f(s) = 0\} = 0$, then $Jf$ is a singleton satisfying*

$$(Jf)(s) = \begin{cases} \|f\|_1, & for\ f(s) > 0, \\ -\|f\|_1, & for\ f(s) < 0, \end{cases} \text{ for all } s \in S. \tag{4.3}$$

(iii) *In particular, if $f \in L_1^{++}(S)$, then $Jf$ is a constant function satisfying*

$$(Jf)(s) = \|f\|_1, \text{ for all } s \in S.$$

*Proof.* Part (i) follows from the representation (4.2) immediately. We show (ii). Assume that $f$ satisfies $\mu\{s \in S: f(s) = 0\} = 0$. Then, $\mu\{s \in S: f(s) > 0\} > 0$, or, $\mu\{s \in S: f(s) < 0\} > 0$, or both. Suppose $\mu\{s \in S: f(s) > 0\} > 0$ (we can similarly prove the case if $\mu\{s \in S: f(s) < 0\} > 0$). Assume that there is $h^* \in L_\infty(S)$ such that $h^*$ is also a value of the normalized duality mapping at $f \in L_1(S)$ that is different from $Jf$ defined by (4.3). Then, $\|h^*\|_\infty = \|f\|_1$, and, for $\mu$-almost $s \in S$, $|h^*(s)| \leq \|h^*\|_\infty = \|f\|_1$. Then, by the assumption that $h^* \neq Jf$ defined by (4.3), there are two possible cases that at least one of them happens.

Case 1. There is $D \in \mathcal{A}$ and $D \subseteq \{s \in S: f(s) > 0\}$ with $\mu(D) > 0$ such that

$$-\|f\|_1 \leq h^*(s) < \|f\|_1, \text{ for all } s \in D.$$

We have

$$\langle h^*, f \rangle = \int_S f(s) h^*(s) \mu(ds)$$

$$= \int_{S \setminus D} f(s) h^*(s) \mu(ds) + \int_D f(s) h^*(s) \mu(ds)$$

$$< \int_{S \setminus D} f(s) h^*(s) \mu(ds) + \int_D f(s) \|f\|_1 \mu(ds)$$

$$\leq \int_{S \setminus D} |f(s) h^*(s)| \mu(ds) + \int_D f(s) \|f\|_1 \mu(ds)$$

$$\leq \int_{S \setminus D} |f(s)| \|f\|_1 \mu(ds) + \int_D f(s) \|f\|_1 \mu(ds)$$

$$= \int_{S \setminus D} |f(s)| \|f\|_1 \mu(ds) + \int_D |f(s)| \|f\|_1 \mu(ds)$$

$$= \|f\|_1 \int_S |f(s)| \mu(ds)$$

$$= \|f\|_1^2.$$

This contradicts to the assumption that $h^*$ is a valued of the normalized duality mapping at $f$.

Case 1. There is $E \in \mathcal{A}$ and $E \subseteq \{s \in S: f(s) < 0\}$ with $\mu(E) > 0$ such that

$$-\|f\|_1 < h^*(s) \le \|f\|_1, \text{ for all } s \in E.$$

We have

$$\langle h^*, f \rangle = \int_S f(s) h^*(s) \mu(ds)$$

$$= \int_{S \setminus E} f(s) h^*(s) \mu(ds) + \int_E f(s) h^*(s) \mu(ds)$$

$$< \int_{S \setminus E} f(s) h^*(s) \mu(ds) + \int_E |f(s)| \|f\|_1 \mu(ds)$$

$$\le \int_{S \setminus E} |f(s) h^*(s)| \mu(ds) + \int_E |f(s)| \|f\|_1 \mu(ds)$$

$$= \int_{S \setminus E} |f(s)| \|f\|_1 \mu(ds) + \int_E |f(s)| \|f\|_1 \mu(ds)$$

$$= \|f\|_1 \int_S |f(s)| \mu(ds)$$

$$= \|f\|_1^2.$$

This contradicts to the assumption that $h^*$ is a valued of the normalized duality mapping at $f$. □

Let $rca(S, \mathcal{A})$ be the Banach space of all regular countable additive measures contained in $ba(S, \mathcal{A})$. Then, $rca(S, \mathcal{A})$ is a closed linear subspace of $ba(S, \mathcal{A})$ (see pages 160 and 162 of section III 7 in [8] by Dunford and Schwartz).

For any $f \in L_1^+(S)$, with $f \ne 0$, $f$ induces its corresponding member $f^{**} \in (S, \mathcal{A})$ as follows

$$f^{**}(A) = \int_A f(s) \mu(ds), \text{ for every } A \in \mathcal{A}.$$

Then, $f^{**}$ is a positive and complete measure on $(S, \mathcal{A})$; and therefore, $f^{**} \in rca(S, \mathcal{A})$, which implies $f^{**} \in ba(S, \mathcal{A}) = L_\infty^*(S)$. For any $k^* \in L_\infty(S)$, it satisfies

$$\langle f^{**}, k^* \rangle_* = \langle k^*, f \rangle = \int_S k^*(s) f(s) \mu(ds).$$

The identification of $f^{**}$ and $f$ embeds $L_1^+(S)$ into $rca(S, \mathcal{A}) \subseteq ba(S, \mathcal{A})$. That is,

$$L_1^+(S) \subseteq rca(S, \mathcal{A}) \subseteq ba(S, \mathcal{A}) = L_\infty^*(S).$$

The results of the following lemma may be known. We provide a simple proof here.

**Lemma 4.3.** $L_1(S)$ *is not strictly convex.*

*Proof.* Take arbitrarily $A, B \in \mathcal{A}$ satisfying $\mu(A) > 0$, $\mu(B) > 0$ and $A \cap B = \emptyset$. Let $\chi_A$ and $\chi_B$ be the indicator functions of the subsets $A$ and $B$, respectively. Define $f = \frac{\chi_A}{\mu(A)}$ and $g = \frac{\chi_B}{\mu(B)}$. Then $\|f\|_1 = \|g\|_1 = 1$. We can check that $\left\|\frac{f+g}{2}\right\|_1 = 1$. □

**Lemma 4.4.** *The normalized duality mapping $J^*: L_\infty(S) \to L_\infty^*(S)$ is not a single valued mapping.*

*Proof.* Assume, by the way of contradiction, that the normalized duality mapping $J^*: L_\infty(S) \to L_\infty^*(S)$ is a single valued mapping. By Theorem 4.3.2 in [28] by Takahashi, $L_\infty(S)$ is smooth. Since $L_\infty(S) = L_1^*(S)$, by Problem 2 in section 4.3 in [28], it implies that $L_1(S)$ is strictly convex. This contradicts to the results of Lemma 4.3. □

We partially showed the well-known properties that all $L_1(S)$, $L_\infty(S)$ and $ba(S, \mathcal{A})$ are neither smooth, nor strictly convex. Where, the proof of the result that $ba(S, \mathcal{A})$ is neither smooth, nor strictly convex is by the properties of $L_\infty(S)$ and Problems 1 and 2 in section 4.3 in [28].

By Lemma 4.1, the normalized duality mapping $J: L_1(S) \rightrightarrows L_\infty(S) = L_1^*(S)$ is not a single valued mapping. Recall that $L_1^*(S) = L_\infty(S)$ and $L_\infty^*(S) = ba(S, \mathcal{A})$. Let $f \in L_1(S)$. Take $f^* \in J(f) \subseteq L_1^*(S) = L_\infty(S)$, by (1.2), the Mordukhovich derivative of $J$ at $(f, f^*)$ is a set valued mapping $\widehat{D}^*J(f, f^*): L_\infty^*(S) \rightrightarrows L_1^*(S)$, which can be rewritten as $\widehat{D}^*J(f, f^*): ba(S, \mathcal{A}) \rightrightarrows L_\infty(S)$. For any $\gamma \in ba(S, \mathcal{A})$, we have

$$\widehat{D}^*J(f,f^*)(\gamma) = \left\{ k^* \in L_1^*(S): \limsup_{\substack{(g,g^*) \to (f,f^*) \\ g \in L_1(S) \text{ and } g^* \in J(g)}} \frac{\langle k^*, g-f \rangle - \langle \gamma, g^*-f^* \rangle_*}{\|g-f\|_1 + \|g^*-f^*\|_\infty} \leq 0 \right\}$$

$$= \left\{ k^* \in L_\infty(S): \limsup_{\substack{(g,g^*) \to (f,f^*) \\ g \in L_1(S) \text{ and } g^* \in J(g)}} \frac{\langle k^*, g-f \rangle - \langle \gamma, g^*-f^* \rangle_*}{\|g-f\|_1 + \|g^*-f^*\|_\infty} \leq 0 \right\}. \quad (4.4)$$

In particular, by $L_1^+(S) \subseteq ba(S, \mathcal{A})$, for any $h \in L_1^+(S)$ that induces $h^{**} \in ba(S, \mathcal{A}) = L_\infty^*(S)$, by (4.4), we have,

$$\widehat{D}^*J(f,f^*)(h^{**}) = \left\{ k^* \in L_\infty(S): \limsup_{\substack{(g,g^*) \to (f,f^*) \\ g \in L_1(S) \text{ and } g^* \in J(g)}} \frac{\langle k^*, g-f \rangle - \langle h^{**}, g^*-f^* \rangle_*}{\|g-f\|_1 + \|g^*-f^*\|_\infty} \leq 0 \right\}$$

$$= \left\{ k^* \in L_\infty(S): \limsup_{\substack{(g,g^*) \to (f,f^*) \\ g \in L_1(S) \text{ and } g^* \in J(g)}} \frac{\langle k^*, g-f \rangle - \langle g^*-f^*, h \rangle}{\|g-f\|_1 + \|g^*-f^*\|_\infty} \leq 0 \right\}. \quad (4.5)$$

**Theorem 4.5.** *Let $f \in L_1(S)$ with $\mu\{s \in S: f(s) = 0\} = 0$. Then, $J(f)$ is a singleton and satisfies*

$$\widehat{D}^*J(f, J(f))(\theta^{**}) = \{\theta^*\}.$$

*Proof.* For any $f \in L_1(S)$, if $\mu\{s \in S: f(s) = 0\} = 0$, then, by part (ii) of Corollary 4.2, $J(f)$ is a singleton, which is denoted by $f^* = J(f)$. It is clear to see that

$$\theta^* \in \widehat{D}^*J(f, f^*)(\theta^{**}). \tag{4.6}$$

Next, we prove that, for any $k^* \in L_\infty(S)$, if $k^* \neq \theta^*$, then

$$k^* \notin \widehat{D}^*J(f, f^*)(\theta^{**}). \tag{4.7}$$

The proof of (4.7) is divided into two cases.

Case 1. $\langle k^*, f \rangle \neq 0$. We may assume $\langle k^*, f \rangle > 0$. In this case, in the limit in (4.5), we take a special direction $g_t = (1 + t)f$, for $t > 0$ with $t \downarrow 0$. By property ($J_4$) of the normalized duality mapping and $f^* = J(f)$, we have

$$(1+t)f^* \in J((1+t)f) = J(g_t), \text{ for all } t > 0.$$

Then, we write

$$g_t^* = (1+t)f^* \in J\big((1+t)f\big) = J(g_t), \text{ for all } t > 0.$$

This implies

$$(g_t, g_t^*) = (g_t, (1+t)f^*) \to (f, f^*), \text{ as } t \downarrow 0.$$

We have

$$\limsup_{\substack{(g,g^*) \to (f,f^*) \\ g \in L_1(S) \text{ and } g^* \in J(g)}} \frac{\langle k^*, g-f \rangle - \langle g^*-f^*, \theta^{**} \rangle_*}{\|g-f\|_1 + \|g^*-f^*\|_\infty}$$

$$= \limsup_{\substack{(g,g^*) \to (f,f^*) \\ g \in L_1(S) \text{ and } g^* \in J(g)}} \frac{\langle k^*, g-f \rangle}{\|g-f\|_1 + \|g^*-f^*\|_\infty}$$

$$\geq \limsup_{\substack{(g_t, g_t^*) \to (f, f^*) \\ g_t \in L_1(S) \text{ and } g_t^* \in J(g_t)}} \frac{\langle k^*, g_t - f \rangle}{\|g_t - f\|_1 + \|g_t^* - f^*\|_\infty}$$

$$= \limsup_{t \downarrow 0} \frac{\langle k^*, (1+t)f - f \rangle}{\|(1+t)f - f\|_1 + \|(1+t)f^* - f^*\|_\infty}$$

$$= \limsup_{t \downarrow 0} \frac{t\langle k^*, f \rangle}{t\|f\|_1 + t\|f^*\|_\infty}$$

$$= \frac{\langle k^*, f \rangle}{2\|f\|_1}$$

$$> 0.$$

This implies

$$k^* \notin \widehat{D}^*J(f,f^*)(\theta^{**}), \text{ for } \langle k^*, f\rangle > 0. \tag{4.7}$$

In case if, $\langle k^*, f\rangle < 0$, then, in the above proof, we take $g_t = (1-t)f$, for $0 < t < 1$ with $t \downarrow 0$. We can prove

$$k^* \notin \widehat{D}^*J(f,f^*)(\theta^{**}), \text{ for } \langle k^*, f\rangle < 0. \tag{4.7}$$

Case 2. $\langle k^*, f\rangle = 0$. Let $A = \{s \in S: f(s) > 0\}$, $B = \{s \in S: f(s) < 0\}$ and $C = \{s \in S: f(s) = 0\}$. Since $\mu(C) = 0$, which implies $f \neq \theta$ and $\mu(A \cup B) \neq 0$. This implies

$$\mu(\{s \in A: k^*(s) > 0\} \cup \{s \in A: k^*(s) < 0\} \cup \{s \in B: k^*(s) > 0\} \cup \{s \in B: k^*(s) < 0\}) \neq 0.$$

At first, we suppose $\mu\{s \in A: k^*(s) < 0\} > 0$. That is, $\mu\{s \in S: f(s) > 0 \text{ and } k^*(s) < 0\} > 0$. Then, there is a positive number $a$ with $a < \|f\|_1$ such that

$$\mu\{s \in S: f(s) > a \text{ and } k^*(s) < 0\} > 0.$$

Let $D = \{s \in S: f(s) > a \text{ and } k^*(s) < 0\}$. By $f \in L_1(S)$, this implies that

$$0 < \mu(D) = \mu\{s \in S: f(s) > a \text{ and } k^*(s) < 0\} < \infty.$$

Let $\chi_D$ be the indicator function of this subset $D$. The above inequality implies that $\chi_D \in L_1(S)$ with $\|\chi_D\|_1 = \mu(D)$. Define

$$h_t = -t\chi_D + f, \text{ for any } t \text{ with } 0 < t < a \text{ and } t \downarrow 0.$$

More precisely, for any $t$ with $0 < t < a$, we have

$$h_t(s) = \begin{cases} f(s) - t, & \text{for } s \in D, \\ f(s), & \text{for } s \notin D. \end{cases} \tag{4.8}$$

By $0 < t < a$ and $f(s) > a$, for all $s \in D$, we have that if $0 < t < a$, then

$$0 < h_t(s) = f(s) - t < f(s), \text{ for all } s \in D.$$

Hence, if $0 < t < a$, then $\|h_t - f\|_1 = t\mu(D)$. This implies that

$$h_t \to f, \text{ as } t \downarrow 0. \tag{4.9}$$

Notice that $a < \|f\|_1$. If $0 < t < a$, we have

$$\|h_t\|_1 = \int_S |h_t(s)|\mu(ds)$$

$$= \int_{S\setminus D}|h_t(s)|\mu(ds) + \int_D |h_t(s)|\mu(ds)$$

$$= \int_{S\setminus D}|f(s)|\mu(ds) + \int_D (f(s) - t)\mu(ds)$$

$$= \int_{S\setminus D}|f(s)|\mu(ds) + \int_D f(s)\mu(ds) - t\mu(D)$$

$$= \int_S |f(s)| \mu(ds) - t\mu(D)$$

$$= \|f\|_1 - t\mu(D).$$

Here, we have $\|f\|_1 - t\mu(D) > 0$. It is because that if if $0 < t < a$, then

$$t\mu(D) < \int_D f(s)\mu(ds) \leq \|f\|_1.$$

It is clear that, by $\mu\{s \in S: f(s) = 0\} = 0$, if $0 < t < a$, then $\mu\{s \in S: h_t(s) = 0\} = 0$. Notice that $D \subseteq A$ and by (4.8), we have $h_t(s) > 0$, for $s \in D \subseteq A$. By Corollary 4.2, we have that $J(h_t)$ is a singleton, which is denoted by $h_t^* = J(h_t)$, for all $0 < t < a$. It satisfies

$$h_t^*(s) = J(h_t)(s) = \begin{cases} \|h_t\|_1, & \text{for } s \in A, \\ -\|h_t\|_1, & \text{for } s \in B \end{cases} = \begin{cases} \|f\|_1 - t\mu(D), & \text{for } s \in A, \\ -(\|f\|_1 - t\mu(D)), & \text{for } s \in B. \end{cases}$$

This implies

$$\|h_t^* - f^*\|_\infty = t\mu(D) \to 0, \text{ as } t \downarrow 0. \tag{4.10}$$

By (4.9) and (4.10), we have

$$(h_t, h_t^*) \to (f, f^*), \text{ as } t \downarrow 0.$$

This implies

$$\limsup_{\substack{(g,g^*) \to (f,f^*) \\ g \in L_1(S) \text{ and } g^* \in J(g)}} \frac{\langle k^*, g-f \rangle - \langle \theta^{**}, g^*-f^* \rangle_*}{\|g-f\|_1 + \|g^*-f^*\|_\infty}$$

$$= \limsup_{\substack{(g,g^*) \to (f,f^*) \\ g \in L_1(S) \text{ and } g^* \in J(g)}} \frac{\langle k^*, g-f \rangle}{\|g-f\|_1 + \|g^*-f^*\|_\infty}$$

$$\geq \limsup_{(h_t,h_t^*) \to (f,f^*)} \frac{\langle k^*, h_t-f \rangle}{\|h_t-f\|_1 + \|h_t^*-f^*\|_\infty}$$

$$= \limsup_{t \downarrow 0} \frac{\langle k^*, -t\chi_D + f - f \rangle}{\|-t\chi_D + f - f\|_1 + \|h_t^*-f^*\|_\infty}$$

$$= \limsup_{t \downarrow 0} \frac{-t\langle k^*, \chi_D \rangle}{t\mu(D) + t\mu(D)}$$

$$= \frac{-\langle k^*, \chi_D \rangle}{2\mu(D)}$$

$$> 0.$$

This implies

$$k^* \notin \widehat{D}^*J(f, J(f))(\theta^{**}), \text{ for } \langle k^*, f \rangle \neq 0 \text{ and } \mu\{s \in A: k^*(s) < 0\} > 0.$$

We can similarly prove that $k^* \notin \widehat{D}^*J(f,J(f))(\theta^{**})$, for $\langle k^*, f \rangle \neq 0$ with respect to any one of the following cases:

$$\mu\{s \in A\colon k^*(s) > 0\} > 0,\ \mu\{s \in B\colon k^*(s) < 0\} > 0,\ \mu\{s \in B\colon k^*(s) > 0\} > 0.$$

This proves (4.7). By (4.6), this theorem is proved. □

**Theorem 4.6.** *For $\theta^* = J(\theta)$, we have*

$$\widehat{D}^*J(\theta, \theta^*)(\theta^{**}) = \{\theta^*\}.$$

*Proof.* It is clear to see that

$$\theta^* \in \widehat{D}^*J(\theta, \theta^*)(\theta^{**}). \tag{4.11}$$

Next, we prove that, for any $k^* \in L_\infty(S)$, if $k^* \neq \theta^*$, then

$$k^* \notin \widehat{D}^*J(\theta, \theta^*)(\theta^{**}). \tag{4.12}$$

By $k^* \neq \theta^*$, we have $\mu(\{s \in S\colon k^*(s) > 0\}) > 0$, or, $\mu\{s \in S\colon k^*(s) < 0\}$, or both. At first, we suppose $\mu\{s \in S\colon k^*(s) > 0\} > 0$. Then, we take $D \subseteq \{s \in S\colon k^*(s) > 0\}$ satisfying

$$0 < \mu(D) < \infty.$$

It follows that $\chi_D \in L_1(S)$ with $\|\chi_D\|_1 = \mu(D)$. Define

$$h_t = t\chi_D, \text{ for any } t > 0 \text{ with } t \downarrow 0.$$

More precisely, for any $t > 0$, we have

$$h_t(s) = \begin{cases} t, & \text{for } s \in D, \\ 0, & \text{for } s \notin D. \end{cases}$$

This implies that, for any $t > 0$, $\|h_t - \theta\|_1 = \|h_t\|_1 = t\mu(D)$. We obtain

$$h_t \to \theta, \text{ as } t \downarrow 0. \tag{4.13}$$

By Lemma 4.1 and $\|h_t\|_1 = t\mu(D)$, we define $j(h_t) \in J(h_t)$ satisfying

$$j(h_t)(s) = \begin{cases} \|h_t\|_1, & \text{for } s \in D, \\ -\|h_t\|_1, & \text{for } s \notin D \end{cases} = \begin{cases} t\mu(D), & \text{for } s \in D, \\ -t\mu(D), & \text{for } s \notin D. \end{cases}$$

This implies

$$\|j(h_t) - \theta^*\|_\infty = \|h_t\|_1 = t\mu(D) \to 0, \text{ as } t \downarrow 0.$$

This is equivalent to

$$j(h_t) \to \theta^*, \text{ as } t \downarrow 0. \tag{4.14}$$

By (4.13) and (4.14), we have

$$(h_t, j(h_t)) \to (\theta, \theta^*), \text{ as } t \downarrow 0.$$

This implies

$$\limsup_{\substack{(g,g^*) \to (\theta,\theta^*) \\ g \in L_1(S) \text{ and } g^* \in J(g)}} \frac{\langle k^*, g-\theta \rangle - \langle \theta^{**}, g^*-\theta^* \rangle_*}{\|g-f\|_1 + \|g^*-\theta^*\|_\infty}$$

$$= \limsup_{\substack{(g,g^*) \to (\theta,\theta^*) \\ g \in L_1(S) \text{ and } g^* \in J(g)}} \frac{\langle k^*, g \rangle}{\|g-f\|_1 + \|g^*\|_\infty}$$

$$\geq \limsup_{(h_t, h_t^*) \to (\theta,\theta^*)} \frac{\langle k^*, h_t \rangle}{\|h_t\|_1 + \|h_t^*\|_\infty}$$

$$= \limsup_{t \downarrow 0} \frac{\langle k^*, t\chi_D \rangle}{\|t\chi_D\|_1 + \|j(h_t)\|_\infty}$$

$$= \limsup_{t \downarrow 0} \frac{t \langle k^*, \chi_D \rangle}{t\mu(D) + t\mu(D)}$$

$$= \frac{\langle k^*, \chi_D \rangle}{2\mu(D)}$$

$$> 0.$$

This proves (4.12). By (4.11) and (4.12), this theorem is proved if $\mu\{s \in S: k^*(s) > 0\} > 0$. We can similarly prove the case if $\mu\{s \in S: k^*(s) < 0\} > 0$. □

Recall that $L_1^+(S)$ is embedded into $rca(S, \mathcal{A}) \subseteq ba(S, \mathcal{A})$. That is, $L_1^+(S) \subseteq rca(S, \mathcal{A}) \subseteq ba(S, \mathcal{A}) = L_\infty^*(S)$. Then, for any $f \in L_1^+(S)$, $f^{**}$ is considered as an element in $ba(S, \mathcal{A})$.

**Theorem 4.7**. *Let* $f \in L_1^+(S)$ *with* $f \neq \theta$ *that induces* $f^{**} \in ba(S, \mathcal{A})$. *Let* $f^* \in J(f)$ *defined by*

$$f^*(s) = \begin{cases} \|f\|_1, & \text{for } f(s) > 0, \\ 0, & \text{for } f(s) = 0, \end{cases} \text{ for all } s \in S. \tag{4.15}$$

*Then, we have*

$$-f^* \notin \widehat{D}^*J(f, f^*)(f^{**}).$$

*Proof.* Let $A = \{s \in S: f(s) > 0\}$, $B = \{s \in S: f(s) < 0\}$ and $C = \{s \in S: f(s) = 0\}$. By the assumption that $f \in L_1^+(S)$ and $f \neq \theta$, we have $\mu(A) > 0$ and $\mu(B) = 0$. Then, there is a positive number $a$ such that

$$\mu\{s \in A: f(s) > a\} > 0.$$

Hence, there is $D \subseteq \{s \in A: f(s) > a\}$ such that $0 < \mu(D) < \infty$. Define

$$h_t = -t\chi_D + f, \text{ for any } 0 < t < a \text{ and } t \downarrow 0.$$

More precisely, for any $t$ with $0 < t < a$, we have

$$h_t(s) = \begin{cases} f(s) - t, & \text{for } s \in D, \\ f(s), & \text{for } s \notin D. \end{cases} \tag{4.16}$$

This implies $\|h_t - f\|_1 = t\mu(D)$, for all $0 < t < a$, and

$$h_t \to f, \text{ as } t \downarrow 0. \tag{4.17}$$

We calculate

$$\|h_t\|_1 = \int_S |h_t(s)| \mu(ds)$$

$$= \int_{S \setminus D} |h_t(s)| \mu(ds) + \int_D |h_t(s)| \mu(ds)$$

$$= \int_{S \setminus D} |f(s)| \mu(ds) + \int_D (f(s) - t) \mu(ds)$$

$$= \int_{S \setminus D} |f(s)| \mu(ds) + \int_D f(s) \mu(ds) - t\mu(D)$$

$$= \int_S |f(s)| \mu(ds) - t\mu(D)$$

$$= \|f\|_1 - t\mu(D). \tag{4.18}$$

Notice that $D \subseteq A$ and by (4.16), we have $h_t(s) > 0$, for $s \in D \subseteq A$. By (4.18), let $h_t^* \in J(h_t)$ be defined as follows (Notice that, for $0 < t < a$, $h_t(s) > 0$ if and only if, $f(s) > 0$, for all $s \in S$).

$$h_t^*(s) = \begin{cases} \|h_t\|_1, & \text{for } h_t(s) > 0, \\ 0, & \text{for } h_t(s) = 0, \end{cases}$$

$$= \begin{cases} \|f\|_1 - t\mu(D), & \text{for } f(s) > 0, \\ 0, & \text{for } f(s) = 0, \end{cases}$$

$$= \begin{cases} \|f\|_1 - t\mu(D), & \text{for } s \in A, \\ 0, & \text{for } s \notin A. \end{cases} \tag{4.19}$$

By (4.15), this implies

$$\|h_t^* - f^*\|_\infty = t\mu(D) \to 0, \text{ as } t \downarrow 0. \tag{4.20}$$

By (4.17) and (4.20), we have

$$(h_t, h_t^*) \to (f, f^*), \text{ as } t \downarrow 0.$$

By (4.15) and the definition of $D$, it is clear to see that

$$\langle f^*, h_t \rangle = \langle f^*, -t\chi_D + f \rangle = \|f\|_1^2 - t\|f\|_1 \mu(D). \tag{4.21}$$

By (4.18) and (4.19), we calculate

$$\langle h_t^*, f \rangle = \int_S h_t^*(s) f(s) \mu(ds)$$

$$= \int_A h_t^*(s)f(s)\,\mu(ds) + \int_{S\setminus A} h_t^*(s)f(s)\,\mu(ds)$$

$$= \int_A (\|f\|_1 - t\mu(D))|f(s)|\,\mu(ds)$$

$$= (\|f\|_1 - t\mu(D))\int_S |f(s)|\,\mu(ds)$$

$$= \|f\|_1^2 - t\|f\|_1\mu(D). \tag{4.22}$$

By (4.21) and (4.22), we have

$$\limsup_{\substack{(g,g^*)\to(f,f^*)\\ g\in L_1(S) \text{ and } g^*\in J(g)}} \frac{\langle -f^*, g-f\rangle - \langle f^{**}, g^*-f^*\rangle_*}{\|g-f\|_1 + \|g^*-f^*\|_\infty}$$

$$= \limsup_{\substack{(g,g^*)\to(f,f^*)\\ g\in L_1(S) \text{ and } g^*\in J(g)}} \frac{-\langle f^*, g-f\rangle - \langle g^*-f^*, f\rangle}{\|g-f\|_1 + \|g^*-f^*\|_\infty}$$

$$\geq \limsup_{\substack{(h_t, h_t^*)\to(f,f^*)\\ h_t\in L_1(S) \text{ and } h_t^*\in J(h_t)}} \frac{-\langle f^*, h_t-f\rangle - \langle h_t^*-f^*, f\rangle}{\|h_t-f\|_1 + \|h_t^*-f^*\|_\infty}$$

$$= \limsup_{t\downarrow 0} \frac{-\langle f^*, h_t\rangle - \langle h_t^*, f\rangle + 2\langle f^*, f\rangle}{\|t\chi_D + f - f\|_1 + \|h_t^*-f^*\|_\infty}$$

$$= \limsup_{t\downarrow 0} \frac{-\langle f^*, h_t\rangle - \langle h_t^*, f\rangle + 2\|f\|_1^2}{t\mu(D) + t\mu(D)}$$

$$= \limsup_{t\downarrow 0} \frac{-\langle f^*, -t\chi_D + f\rangle - \langle h_t^*, f\rangle + 2\|f\|_1^2}{2t\mu(D)}$$

$$= \limsup_{t\downarrow 0} \frac{-(\|f\|_1^2 - t\|f\|_1\mu(D)) - (\|f\|_1^2 - t\|f\|_1\mu(D)) + 2\|f\|_1^2}{2t\mu(D)}$$

$$= \limsup_{t\downarrow 0} \frac{2t\|f\|_1\mu(D)}{2t\mu(D)}$$

$$= \|f\|_1$$

$$> 0.$$

This implies that $-f^* \notin \widehat{D}^*J(f, f^*)(f^{**})$. This theorem is proved. □

We define an ordering relation $\prec$ on $L_\infty(S)$ as follows. For any $f^*, u^* \in L_\infty(S)$, we write

$$f^* \prec u^* \text{ if and only if } f^*(s) < u^*(s), \text{ for almost all } s \in S.$$

For any $f^* \in L_\infty(S)$, we write

$$(f^*)_\prec = \{u^* \in L_\infty(S): f^* \prec u^*\}.$$

**Corollary 4.8.** *Let $f \in L_1^{++}(S)$ that induces $f^{**} \in ba(S, \mathcal{A})$. Then $J(f)$ is a singleton satisfying*

$$\widehat{D}^*J(f,J(f))(f^{**}) \cap (J(f))_{\prec} = \emptyset.$$

*Proof.* Since $f \in L_1^{++}(S)$, let $J(f) = f^*$, which is a constant function with value $\|f\|_1$. For any $u^* \in L_\infty(S)$ with $f^* \prec u^*$, there is a positive number $b$ such that, there is $E \subseteq \{s \in S : u^*(s) > \|f\|_1 + b\}$ satisfying $0 < \mu(E) < \infty$. Define

$$h_t = t\chi_E + f, \text{ for any } t > 0 \text{ with } t \downarrow 0.$$

More precisely, for any $t > 0$, we have

$$h_t(s) = \begin{cases} f(s) + t, & \text{for } s \in E, \\ f(s), & \text{for } s \notin E. \end{cases}$$

This implies $\|h_t - f\|_1 = t\mu(E)$, for all $t > 0$, and $h_t \to f$, as $t \downarrow 0$. Similar to (4.8), we calculate

$$\|h_t\|_1 = \int_S |h_t(s)|\mu(ds)$$

$$= \int_{S\setminus E} |f(s)|\mu(ds) + \int_E (f(s) + t)\mu(ds)$$

$$= \int_{S\setminus E} |f(s)|\mu(ds) + \int_E f(s)\mu(ds) + t\mu(E)$$

$$= \int_S |f(s)|\mu(ds) + t\mu(E)$$

$$= \|f\|_1 + t\mu(E).$$

By $f \in L_1^{++}(S)$, we have $h_t(s) \in L_1^{++}(S)$, for all $t > 0$. Then $h_t^* := J(h_t)$ is a singleton with value of $\|h_t\|_1$. That is,

$$h_t^* := J(h_t) = \|h_t\|_1 = \|f\|_1 + t\mu(E), \text{ for all } t > 0.$$

This implies $\|h_t^* - f^*\|_\infty = t\mu(E) \to 0$, as $t \downarrow 0$. We have

$$(h_t, h_t^*) \to (f, f^*), \text{ as } t \downarrow 0.$$

It is clear to see that $\langle f^*, h_t \rangle = \langle f^*, t\chi_E + f \rangle = \|f\|_1^2 + t\|f\|_1\mu(E)$. We calculate

$$\langle h_t^*, f \rangle = \int_S h_t^*(s)f(s)\mu(ds)$$

$$= (\|f\|_1 + t\mu(E))\int_S |f(s)|\mu(ds)$$

$$= \|f\|_1^2 + t\|f\|_1\mu(E).$$

and $\quad \langle u^*, h_t - f \rangle = \langle u^*, t\chi_E \rangle$

$$= t \int_E u^*(s)\,\mu(ds)$$

$$> t \int_E (f^*(s) + b)\,\mu(ds)$$

$$= t(\|f\|_1 + b)\mu(E).$$

We have

$$\limsup_{\substack{(g,g^*) \to (f,f^*) \\ g \in L_1(S) \text{ and } g^* \in J(g)}} \frac{\langle u^*, g-f \rangle - \langle f^{**},\ g^*-f^* \rangle_*}{\|g-f\|_1 + \|g^*-f^*\|_\infty}$$

$$= \limsup_{\substack{(g,g^*) \to (f,f^*) \\ g \in L_1(S) \text{ and } g^* \in J(g)}} \frac{\langle u^*, g-f \rangle - \langle g^*-f^*, f \rangle}{\|g-f\|_1 + \|g^*-f^*\|_\infty}$$

$$\geq \limsup_{\substack{(h_t, h_t^*) \to (f,f^*) \\ h_t \in L_1(S) \text{ and } h_t^* \in J(h_t)}} \frac{\langle u^*, h_t-f \rangle - \langle h_t^*-f^*, f \rangle}{\|h_t-f\|_1 + \|h_t^*-f^*\|_\infty}$$

$$= \limsup_{t \downarrow 0} \frac{\langle u^*, t\chi_E \rangle - \langle t\mu(E), f \rangle}{2t\mu(E)}$$

$$> \limsup_{t \downarrow 0} \frac{t(\|f\|_1 + b)\mu(E) - t\|f\|_1 \mu(E)}{2t\mu(E)}$$

$$= \limsup_{t \downarrow 0} \frac{tb\mu(E)}{2t\mu(E)}$$

$$= \frac{b}{2} > 0.$$

This implies that

$$u^* \notin \widehat{D}^*J(f, J(f))(f^{**}), \text{ for any } u^* \in L_\infty(S) \text{ with } J(f) \prec u^*. \qquad \square$$

### 5. The normalized duality mapping in (general) Banach space $C[0, 1]$

Let $(C[0, 1], \|\cdot\|, \Sigma)$ be the Banach space of all continuous real valued functions on $[0, 1]$ with respect to the standard Borel $\sigma$-field $\Sigma$ and with the maximum norm (see [6, 29] for more details)

$$\|f\| = \max_{0 \leq s \leq 1} |f(s)|, \text{ for any } f \in C[0, 1].$$

The dual space of $C[0, 1]$ is denoted by $C^*[0, 1]$ that is $rca[0, 1]$. Let $\langle \cdot,\ \cdot \rangle$ denote the real canonical pairing between $C^*[0, 1]$ and $C[0, 1]$. In this section, let $m$ denote the standard measure on $[0, 1]$ satisfying $m[a, b] = b - a$, for any $0 \leq a < b \leq 1$. It is clear to see $m \in rca[0, 1]$.

By the Riesz Representation Theorem (see Theorem IV.6.3 in Dunford and Schwartz [8]), for any $\varphi \in C^*[0, 1]$, there is a real valued, regular and countable additive functional $\mu \in rca[0, 1]$, which is defined on the given $\sigma$-field $\Sigma$ on $[0, 1]$, such that

$$\langle \varphi, f \rangle = \int_0^1 f(s)\mu(ds), \text{ for any } f \in C[0, 1]. \tag{5.1}$$

Throughout this section, without any special mention, we shall identify $\varphi$ and $\mu$ in (5.1). We say that $\mu \in C^*[0, 1]$ ((it is $rca[0, 1]$) and (5.1) is rewritten as

$$\langle \mu, f \rangle = \int_0^1 f(s)\mu(ds), \text{ for any } f \in C[0, 1]. \tag{5.1}$$

The norm of $\mu \in C^*[0, 1] = rca[0, 1]$ is denoted by $\|\mu\|_*$ that is defined by

$$\|\mu\|_* := v(\mu, [0, 1]).$$

Where, $v(\mu, [0, 1])$ is the total variation of $\mu$ on $[0, 1]$, which is defined by (see page 160 in [8])

$$v(\mu, [0, 1]) = \sup_{E \in \Sigma} |\mu(E)|.$$

The origin of $C[0, 1]$ is denoted by $\theta$, which is the constant function defined on $[0, 1]$ with value 0. The origin of the dual space $C^*[0, 1]$ is denoted by $\theta^*$, which is also a constant functional on $\Sigma$ with value 0. This is,

$$\theta^*(E) = 0, \text{ for every } E \in \Sigma.$$

The dual space of $rca[0, 1]$ is denoted by $rca^*[0, 1]$. The real canonical pairing between the dual space $rca^*[0, 1]$ $(= C^{**}[0, 1])$ and $rca[0, 1]$ $(= C^*[0, 1])$ is denoted by $\langle \cdot, \cdot \rangle_*$. In this section, we use English letters to name the members in $C[0, 1]$, such as, $f$, $g$, …; lower case Greek letters for members in $C^*[0, 1]$; upper case Greek letters for members in $C^{**}[0, 1]$).

According to Dunford and Schwartz (see (F1) on page 374 in [8]), so far, there is no completely satisfactory representation for $rca^*[0, 1]$, which is the conjugate of the space $rca[0, 1]$. However, we could find a subset of $rca^*[0, 1]$ for us to use in this section.

The positive cone of $C[0, 1]$ is denoted by $C^+[0, 1]$ that is defined by

$$C^+[0, 1] = \{f \in C[0, 1]: f(s) \geq 0, \text{ for all } s \in [0, 1]\}.$$

$C^+[0, 1]$ is a pointed closed and convex cone in $C[0, 1]$. We define the strict positive "cone" as follows

$$C^{++}[0, 1] = \{f \in C[0, 1]: f(s) > 0, \text{ for all } s \in [0, 1]\}.$$

$C^{++}[0, 1]$ is a convex and open subset in $C[0, 1]$. However, $\{\theta\} \cup C^{++}[0, 1]$ is a pointed convex cone in $C[0, 1]$ with vertex $\theta$, which is neither closed, nor open. For any given $f \in C^+[0, 1]$, $f$ defines a member $f^{**} \in C^{**}[0, 1] = rca^*[0, 1]$ as follows

$$\langle f^{**}, \gamma \rangle_* = \langle \gamma, f \rangle = \int_0^1 f(s)\gamma(ds), \text{ for any } \gamma \in rca[0, 1]. \tag{5.2}$$

It satisfies

$$|\langle f^{**}, \gamma \rangle_*| \leq \|f\| \|\gamma\|_*, \text{ for any } \gamma \in rca[0, 1].$$

This implies $\|f^{**}\|_{rca^*} = \|f\|$. By this definition, $C^+[0, 1]$ is embedded into $rca^*[0, 1]$. That is

$$C^+[0, 1] \subseteq rca^*[0, 1] = C^{**}[0, 1].$$

For any $f \in C[0, 1]$, we define

$$M(f) = \{t \in [0, 1]: |f(t)| = \|f\|\}.$$

By the continuity of $f$ on $[0, 1]$, $M(f)$ is a nonempty closed subset of $[0, 1]$, which is called the maximizing set of $f$. The connection between the normalized duality mapping $J: C \to 2^{C^*}\setminus\{\emptyset\}$ and the maximizing sets is provided in [18−19].

**Lemma 5.1 [18].** *For any $f \in C[0, 1]$ with $\|f\| > 0$, then*

(a) $\mu \in J(f) \implies v(\mu, [0, 1]\setminus M(f)) = 0$;
(b) *Let $s_j \in M(f)$, for $j = 1, 2, \ldots, m$, for some positive integer $m$, define $\mu \in rca[0, 1]$ by*

$$\mu(s_j) = \alpha_j f(s_j), \text{ for } j = 1, 2, \ldots, m$$

*and*

$$\mu(E) = 0, \text{ for any } E \subseteq [0, 1]\setminus\{s_j: j = 1, 2, \ldots, m\}, \tag{5.3}$$

*where $\alpha_j > 0$, for $j = 1, 2, \ldots, m$ satisfying $\sum_{j=1}^{m} \alpha_j = 1$. Then, $\mu \in J(f)$.*

(c) *If $M(f)$ is not a singleton, then $J(f)$ is an infinite set.*

The following lemmas provide some properties of the maximizing sets.

**Lemma 5.1.** *Let $f \in C[0, 1]$ with $\|f\| > 0$. Suppose that there is $[a, b] \subseteq [0, 1]$ with $a < b$ satisfying $f(s) = \|f\|$, for all $s \in [a, b]$. Define $\mu \in rca[0, 1]$ as follows, for any $E \in \Sigma$,*

$$\mu(E) = \begin{cases} \frac{\|f\|}{b-a} m(E), & \text{if } E \subseteq [a, b], \\ 0, & \text{if } E \cap [a, b] = \emptyset. \end{cases}$$

*Then, $\mu \in J(f)$.*

*Proof.* The proof of this lemma is similar to the proof of Lemma 5.1 [18]. It is omitted here. □

**Lemma 5.2.** *For any $f \in C[0, 1]$ with $\|f\| > 0$ and for any $t \neq 0$, we have*

$$M(tf) = M(f).$$

*Proof.* The proof of this lemma is straightforward and it is omitted here. □

For $(f, \mu) \in \text{gph}J$, that is, for $f \in C[0, 1]$ and $\mu \in J(f) \subseteq rca[0, 1] = C^*[0,1]$, by definition (1.2), the Mordukhovich derivative of $J$ at point $(f, \mu)$ is a set valued mapping $\widehat{D}^*J(f, \mu): C^{**}[0,1] \rightrightarrows C^*[0,1]$. For any $\Phi \in C^{**}[0,1]$, it is defined by

$$\widehat{D}^*J(f,\mu)(\Phi) = \left\{\lambda \in C^*[0,1]: \limsup_{\substack{(g,\gamma) \to (f,\mu) \\ g \in C[0,1] \text{ and } \gamma \in J(g)}} \frac{\langle \lambda, g-f \rangle - \langle \Phi, \gamma-\mu \rangle_*}{\|g-f\| + \|\gamma-\mu\|_*} \leq 0\right\}. \quad (5.4)$$

We use the above definition to investigate some properties of the Mordukhovich derivatives of the set valued normalized duality mapping $J$ on $C[0, 1]$.

**Theorem 5.3**. *Let $f \in C^+[0, 1]$ with $\|f\| > 0$ that induces $f^{**} \in rca^*[0, 1]$ defined in (5.2). Then, for any $\mu \in J(f)$, we have*

$$\theta^* \notin \widehat{D}^*J(f,\mu)(f^{**}).$$

*Proof.* By (5.4), for any given $\mu \in J(f)$, we have

$$\theta^* \in \widehat{D}^*J(f,\mu)(f^{**}) \iff \limsup_{\substack{(g,\gamma) \to (f,\mu) \\ g \in C[0,1] \text{ and } \gamma \in J(g)}} \frac{\langle \theta^*, g-f \rangle - \langle f^{**}, \gamma-\mu \rangle_*}{\|g-f\| + \|\gamma-\mu\|_*} \leq 0$$

$$\iff \limsup_{\substack{(g,\gamma) \to (f,\mu) \\ g \in C[0,1] \text{ and } \gamma \in J(g)}} \frac{-\langle f^{**}, \gamma-\mu \rangle_*}{\|g-f\| + \|\gamma-\mu\|_*} \leq 0. \quad (5.5)$$

We estimate the limit in (5.5). By the definition of $f^{**}$ in (5.2), we have

$$\limsup_{\substack{(g,\gamma) \to (f,\mu) \\ g \in C[0,1] \text{ and } \gamma \in J(g)}} \frac{-\langle f^{**}, \gamma-\mu \rangle_*}{\|g-f\| + \|\gamma-\mu\|_*} = \limsup_{\substack{(g,\gamma) \to (f,\mu) \\ g \in C[0,1] \text{ and } \gamma \in J(g)}} \frac{-\langle \gamma-\mu, f \rangle}{\|g-f\| + \|\gamma-\mu\|_*}. \quad (5.6)$$

In the limit (5.6), we take a special direction $g_t \in C[0,1]$ and $\gamma_t \in J(g_t)$, for $t \downarrow 0$, as follows

$$g_t = (1-t)f \text{ and } \gamma_t = (1-t)\mu, \text{ for } 1 > t > 0.$$

By the properties of the normalized duality mapping and the assumption that $\mu \in J(f)$, we have $\gamma_t \in J(g_t)$ satisfying

$$\|\mu - \gamma_t\|_* = t\|\mu\|_*, \text{ for all } 1 > t > 0.$$

This implies

$$g_t \to f \text{ in } C[0,1] \text{ and } \gamma_t \to \mu \text{ in } rca[0, 1], \text{ as } t \downarrow 0.$$

Calculating the limit in (5.6), we have

$$\limsup_{\substack{(g,\gamma) \to (f,\mu) \\ g \in C[0,1] \text{ and } \gamma \in J(g)}} \frac{-\langle \gamma-\mu, f \rangle}{\|g-f\| + \|\gamma-\mu\|_*}$$

$$\geq \limsup_{(g_t,\gamma_t) \to (f,\mu)} \frac{-\langle \gamma_t-\mu, f \rangle}{\|g_t-f\| + \|\gamma_t-\mu\|_*}$$

$$= \limsup_{t \downarrow 0} \frac{t \langle \mu, f \rangle}{t\|f\| + t\|\mu\|_*}$$

$$= \frac{\|f\|^2}{2\|f\|}$$

$$= \frac{\|f\|}{2} > 0.$$

By (5.5), this implies $\theta^* \notin \widehat{D}^*J(f,\mu)(f^{**})$. □

**Theorem 5.4.** *Let $f \in C[0, 1]$ and $\lambda \in rca[0, 1]$. If $\langle \lambda, f \rangle \neq 0$, then, for any $\mu \in J(f)$, we have*

$$\lambda \notin \widehat{D}^*J(f,\mu)(\theta^{**}).$$

*Proof.* For any given $f \in C[0, 1]$ and $\lambda \in rca[0, 1]$, the condition $\langle \lambda, f \rangle \neq 0$ implies $f \neq \theta$, that is $\|f\| > 0$. By (5.4), for any $\mu \in J(f)$, we have

$$\lambda \in \widehat{D}^*J(f,\mu)(\theta^{**}) \Leftrightarrow \limsup_{\substack{(g,\gamma) \to (f,\mu) \\ g \in C[0,1] \text{ and } \gamma \in J(g)}} \frac{\langle \lambda, g-f \rangle - \langle \theta^{**}, \gamma - \mu \rangle_*}{\|g-f\| + \|\gamma - \mu\|_*} \leq 0$$

$$\Leftrightarrow \limsup_{\substack{(g,\gamma) \to (f,\mu) \\ g \in C[0,1] \text{ and } \gamma \in J(g)}} \frac{\langle \lambda, g-f \rangle}{\|g-f\| + \|\gamma - \mu\|_*} \leq 0. \tag{5.7}$$

For any $0 < t < 1$, let

$$g_t = \begin{cases} (1+t)f, & \text{if } \langle \lambda, f \rangle > 0, \\ (1-t)f, & \text{if } \langle \lambda, f \rangle < 0. \end{cases}$$

and

$$\gamma_t = \begin{cases} (1+t)\mu, & \text{if } \langle \lambda, f \rangle > 0, \\ (1-t)\mu, & \text{if } \langle \lambda, f \rangle < 0. \end{cases}$$

By $\mu \in J(f)$, we have $\gamma_t \in J(g_t)$, for any $0 < t < 1$. Calculating the limit in (5.7), we have

$$\limsup_{\substack{(g,\gamma) \to (f,\mu) \\ g \in C[0,1] \text{ and } \gamma \in J(g)}} \frac{\langle \lambda, g-f \rangle}{\|g-f\| + \|\gamma - \mu\|_*}$$

$$\geq \limsup_{(g_t,\gamma_t) \to (f,\mu)} \frac{\langle \lambda, g_t-f \rangle}{\|g_t-f\| + \|\gamma_t - \mu\|_*}$$

$$= \limsup_{t \downarrow 0} \frac{t|\langle \lambda, f \rangle|}{t\|f\| + t\|\mu\|_*}$$

$$= \frac{|\langle \lambda, f \rangle|}{2\|f\|} > 0.$$

By (5.7), this implies $\lambda \notin \widehat{D}^*J(f,\mu)(\theta^{**})$. □

**Theorem 5.5.** *Let $\lambda \in rca[0, 1]$. If $\lambda[0, 1] \neq 0$, then for any $f \in C[0, 1]$ there is $\mu \in J(f)$ such that*

$$\lambda \notin \widehat{D}^*J(f,\mu)(\theta^{**}).$$

*In particular, if $\lambda[0, 1] \neq 0$, then*

$$\lambda \notin \widehat{D}^*J(\theta,\theta^*)(\theta^{**}).$$

*Proof.* Let $\lambda \in rca[0, 1]$ with $\lambda[0, 1] \neq 0$. The proof of this theorem is divided into two cases.

Case 1. $\lambda[0, 1] > 0$. The proof of this case is also divided into two subcases.

Subcase 1.1. There is a subset of the maximizing set of $f$, $\{s_j \in M(f): \text{for } j = 1, 2, \ldots, m\}$, for some positive integer $m$ such that

$$f(s_j) = \|f\|, \text{ for } j = 1, 2, \ldots, m. \tag{5.8}$$

In this case, we define $\mu \in rca[0, 1]$ by

$$\mu(s_j) = \alpha_j f(s_j), \text{ for } j = 1, 2, \ldots, m$$

and

$$\mu(E) = 0, \text{ for any } E \subseteq [0, 1] \setminus \{s_j : j = 1, 2, \ldots, m\}. \tag{5.3}$$

where $\alpha_j > 0$, for $j = 1, 2, \ldots, m$ satisfying $\sum_{j=1}^{m} \alpha_j = 1$. By Lemma 5.1 in [8], we have $\mu \in J(f)$. For any given $t > 0$, we define $h_t \in C[0,1]$ by

$$h_t(s) = t + f(s), \text{ for all } s \in [0,1]. \tag{5.9}$$

We have

$$\|h_t\| = t + \|f\| \quad \text{and} \quad \{s_j : \text{for } j = 1, 2, \ldots, m\} \subseteq M(h_t), \text{ for any } t > 0. \tag{5.10}$$

Then, for any $t > 0$, similar to (5.3), we define $\mu_t \in rca[0, 1]$ by

$$\mu_t(s_j) = \alpha_j h_t(s_j) = \alpha_j(t + f(s_j)), \text{ for } j = 1, 2, \ldots, m$$

and

$$\mu_t(E) = 0, \text{ for any } E \subseteq [0, 1] \setminus \{s_j : j = 1, 2, \ldots, m\}, \tag{5.11}$$

where $\alpha_j$ is given in (5.3), for $j = 1, 2, \ldots, m$. Then, $\mu_t \in J(h_t)$. By (5.3) and (5.11), we have

$$\|\mu_t - \mu\|_* = t, \text{ for any } t > 0. \tag{5.12}$$

By (5.9) and (5.12), we obtain

$$(h_t, \mu_t) \to (f, \mu), \text{ as } t \downarrow 0. \tag{5.13}$$

By (5.12) and (5.13), calculating the limit in (5.7), we have

$$\limsup_{\substack{(g,\gamma) \to (f,\mu) \\ g \in C[0,1] \text{ and } \gamma \in J(g)}} \frac{\langle \lambda, g-f \rangle}{\|g-f\| + \|\gamma-\mu\|_*}$$

$$\geq \limsup_{(h_t,\mu_t) \to (f,\mu)} \frac{\langle \lambda, h_t-f \rangle}{\|h_t-f\| + \|\mu_t-\mu\|_*}$$

$$= \limsup_{t \downarrow 0} \frac{t\lambda[0,1]}{2t}$$

$$= \frac{\lambda[0,1]}{2} > 0.$$

By (5.7), this implies $\lambda \notin \widehat{D}^*J(f,\mu)(\theta^{**})$.

Subcase 1.2. Suppose that such a subset satisfying (5.8) of the maximizing set of $f$ does not exist. That is,

$$f(s) < \|f\|, \text{ for all } s \in [0,1]. \tag{5.14}$$

This implies

$$A := \max\{f(s) : s \in [0,1]\} < \|f\|.$$

In this case, there is a subset of the maximizing set of $f$, $\{s_k \in M(f): \text{ for } k = 1, 2, \ldots, n\}$, for some positive integer $n$ such that

$$f(s_k) = -\|f\|, \text{ for } k = 1, 2, \ldots, n.$$

Then, we define $\mu \in rca[0,1]$ by

$$\mu(s_k) = \beta_k f(s_k), \text{ for } k = 1, 2, \ldots, n$$

and

$$\mu(E) = 0, \text{ for any } E \subseteq [0,1] \setminus \{s_k : j = 1, 2, \ldots, n\}. \tag{5.15}$$

where $\beta_k > 0$, for $k = 1, 2, \ldots, n$ satisfying $\sum_{k=1}^{n} \beta_k = 1$. By Lemma 5.1 in [8], we have $\mu \in J(f)$. For any positive number $t$ with $0 < t < \frac{1}{2}(\|f\| - A)$, we define $u_t \in C[0,1]$ by

$$u_t(s) = t + f(s), \text{ for all } s \in [0,1].$$

This implies

$$u_t(s_k) = t + f(s_k) = -\|f\| + t, \text{ for } k = 1, 2, \ldots, n.$$

For any $0 < t < \frac{1}{2}(\|f\| - A)$ and for any $k = 1, 2, \ldots, n$, we have

$$f(s_k) + 2t < -\|f\| + (\|f\| - A) = -A, \text{ for } k = 1, 2, \ldots, n.$$

This implies that, for any $0 < t < \frac{1}{2}(\|f\| - A)$, we have

$$f(s_k) + t < -(A+t), \text{ for } k = 1, 2, \ldots, n.$$

Then, for any $0 < t < \frac{1}{2}(\|f\| - A)$ and for any $k = 1, 2, \ldots, n$, $u_t(s_k)$ satisfies

$$u_t(s_k) = -\|f\| + t = f(s_k) + t \leq f(s) + t = u_t(s), \text{ if } f(s) \leq 0;$$

and

$$u_t(s_k) = -\|f\| + t = f(s_k) + t < -(A + t) \leq -(f(s) + t) = -u_t(s), \text{ if } f(s) > 0.$$

This implies that, for any $0 < t < \frac{1}{2}(\|f\| - A)$,

(a) $\|u_t\| = \|f\| - t$;
(b) $\{s_k : \text{for } k = 1, 2, \ldots, n\} \subseteq M(u_t)$.

By (5.15), for any $0 < t < \frac{1}{2}(\|f\| - A)$, we define $\mu_t \in rca[0, 1]$ by

$$\mu_t(s_k) = \beta_k(f(s_k) + t) = \beta_k(u_t(s_k)), \text{ for } k = 1, 2, \ldots, n$$

and

$$\mu_t(E) = 0, \text{ for any } E \subseteq [0, 1] \setminus \{s_k : k = 1, 2, \ldots, n\}. \tag{5.16}$$

where $\beta_k$ satisfies (5.15). Then, $\mu_t \in J(u_t)$, for any $0 < t < \frac{1}{2}(\|f\| - A)$. We have

$$\|u_t\| = \|f\| - t \quad \text{and} \quad \|u_t - f\| = t, \text{ for any } 0 < t < \frac{1}{2}(\|f\| - A). \tag{5.17}$$

By (5.15) and (5.16), we have

$$\|\mu_t - \mu\|_* = t, \text{ for any } 0 < t < \frac{1}{2}(\|f\| - A). \tag{5.18}$$

By (5.17) and (5.18), we obtain

$$(u_t, \mu_t) \to (f, \mu), \text{ as } t \downarrow 0. \tag{5.19}$$

By (5.18) and (5.19), calculating the limit in (5.7), we have

$$\limsup_{\substack{(g,\gamma) \to (f,\mu) \\ g \in C[0,1] \text{ and } \gamma \in J(g)}} \frac{\langle \lambda, g - f \rangle}{\|g - f\| + \|\gamma - \mu\|_*}$$

$$\geq \limsup_{(u_t, \mu_t) \to (f,\mu)} \frac{\langle \lambda, u_t - f \rangle}{\|u_t - f\| + \|\mu_t - \mu\|_*}$$

$$= \limsup_{t \downarrow 0} \frac{t\lambda[0,1]}{2t}$$

$$= \frac{\lambda[0,1]}{2} > 0.$$

By (5.7), this implies $\lambda \notin \widehat{D}^*J(f, \mu)(\theta^{**})$.

Case 2. $\lambda[0, 1] < 0$. Similar to the proof of case 1, the proof of Case 2 is also divided into two subcases.

Subcase 2.1. There is a subset of the maximizing set of $f$, $\{s_j \in M(f): \text{for } j = 1, 2, \ldots, m\}$, for some positive integer $m$ such that

$$f(s_j) = -\|f\|, \text{ for } j = 1, 2, \ldots, m.$$

In this case, we define $h_t \in C[0,1]$ by

$$h_t(s) = f(s) - t, \text{ for all } s \in [0,1].$$

Then, rest of the proof of subcase 2.1 is similar to the proof of subcase 1.1.

Subcase 2.2. Suppose that

$$f(s) > -\|f\|, \text{ for all } s \in [0,1].$$

This implies

$$B := \min\{f(s): s \in [0,1]\} > -\|f\|.$$

In this case, for any positive number $t$ with $0 < t < \frac{1}{2}(\|f\| + B)$, we define $u_t \in C[0,1]$ by

$$u_t(s) = f(s) - t, \text{ for all } s \in [0,1].$$

In this case, there is a subset of the maximizing set of $f$, $\{s_k \in M(f): \text{for } k = 1, 2, \ldots, n\}$, for some positive integer $n$ such that

$$f(s_k) = \|f\|, \text{ for } k = 1, 2, \ldots, n.$$

This implies

$$u_t(s_k) = f(s_k) - t = \|f\| - t, \text{ for } k = 1, 2, \ldots, n.$$

For any $0 < t < \frac{1}{2}(\|f\| + B)$ and for any $k = 1, 2, \ldots, n$, we have

$$f(s_k) - 2t > -\|f\| - (\|f\| + B) = -B, \text{ for } k = 1, 2, \ldots, n.$$

This implies that, for any $0 < t < \frac{1}{2}(\|f\| + B)$, we have

$$f(s_k) - t > -B + t, \text{ for } k = 1, 2, \ldots, n.$$

Then, for any $0 < t < \frac{1}{2}(\|f\| + B)$ and for any $k = 1, 2, \ldots, n$, $u_t(s_k)$ satisfies

$$u_t(s_k) = \|f\| - t = f(s_k) - t \geq f(s) - t = u_t(s), \text{ if } f(s) \geq 0;$$

and

$$u_t(s_k) = \|f\| - t = f(s_k) - t > -B + t \geq -f(s) + t = -u_t(s), \text{ if } f(s) < 0.$$

This implies that, for any $0 < t < \frac{1}{2}(\|f\| + B)$,

(a) $\|u_t\| = \|f\| - t$;
(b) $\{s_k: \text{for } k = 1, 2, \ldots, n\} \subseteq M(u_t)$.

Then, rest of the proof of subcase 2.2 is similar to the proof of subcase 1.2. □

**Theorem 5.6**. *Let $f \in C^+[0, 1]$ with $\|f\| > 0$ that induces $f^{**} \in rca^*[0, 1]$. For any $u \in C^+[0, 1]$, if $\|u\| > \|f\|$ and $M(u) \cap M(f) \neq \emptyset$, then, there are $\mu \in J(f)$ and $\lambda \in J(u)$ such that*

$$\lambda \notin \widehat{D}^* J(f, \mu)(f^{**}).$$

*Proof.* For any given $u \in C^+[0, 1]$, suppose $\|u\| > \|f\|$ and $M(u) \cap M(f) \neq \emptyset$. Then, there is a set $\{s_j \in M(u) \cap M(f): \text{for } j = 1, 2, \ldots, m\}$, for some positive integer $m$ such that

$$u(s_j) = \|g\| \quad \text{and} \quad f(s_j) = \|f\|, \text{ for } j = 1, 2, \ldots, m. \tag{5.20}$$

We take $\alpha_j > 0$, for $j = 1, 2, \ldots, m$ satisfying $\sum_{j=1}^{m} \alpha_j = 1$. We define $\mu \in rca[0, 1]$ by

$$\mu(s_j) = \alpha_j f(s_j), \text{ for } j = 1, 2, \ldots, m$$

and

$$\mu(E) = 0, \text{ for any } E \subseteq [0, 1] \setminus \{s_j : j = 1, 2, \ldots, m\}.$$

And, we define $\lambda \in rca[0, 1]$ by

$$\lambda(s_j) = \alpha_j u(s_j), \text{ for } j = 1, 2, \ldots, m$$

and

$$\lambda(E) = 0, \text{ for any } E \subseteq [0, 1] \setminus \{s_j : j = 1, 2, \ldots, m\}.$$

By Lemma 5.1 in [8], we have $\mu \in J(f)$ and $\lambda \in J(u)$. For any given $t > 0$, we define $h_t$ by

$$h_t(s) = t + f(s), \text{ for all } s \in [0,1].$$

Similar to (5.20), for any given $t > 0$, we define $\mu_t \in rca[0, 1]$ by

$$\mu_t(s_j) = \beta_j(f(s_j) + t), \text{ for } j = 1, 2, \ldots, m$$

and

$$\mu_t(E) = 0, \text{ for any } E \subseteq [0, 1] \setminus \{s_j : j = 1, 2, \ldots, m\}.$$

Then, $\mu_t \in J(u_t)$, for any $t > 0$. Similar to the proof of the previous theorems, we have

$$(h_t, \mu_t) \to (f, \mu), \text{ as } t \downarrow 0. \tag{5.21}$$

By (5.7), we have

$$\lambda \in \widehat{D}^* J(f, \mu)(f^{**}) \Leftrightarrow \limsup_{\substack{(g,\gamma) \to (f,\mu) \\ g \in C[0,1] \text{ and } \gamma \in J(g)}} \frac{\langle \lambda, g-f \rangle - \langle f^{**}, \gamma-\mu \rangle_*}{\|g-f\| + \|\gamma-\mu\|_*} \leq 0$$

$$\Leftrightarrow \limsup_{\substack{(g,\gamma) \to (f,\mu) \\ g \in C[0,1] \text{ and } \gamma \in J(g)}} \frac{\langle \lambda, g-f \rangle - \langle \gamma-\mu, f \rangle}{\|g-f\| + \|\gamma-\mu\|_*} \leq 0. \tag{5.22}$$

By (5.21), we calculate the limit in (5.22).

$$\limsup_{\substack{(g,\gamma) \to (f,\mu) \\ g \in C[0,1] \text{ and } \gamma \in J(g)}} \frac{\langle \lambda, g-f \rangle - \langle \gamma-\mu, f \rangle}{\|g-f\| + \|\gamma-\mu\|_*}$$

$$\geq \limsup_{(h_t,\mu_t) \to (f,\mu)} \frac{\langle \lambda, h_t-f \rangle - \langle \mu_t-\mu, f \rangle}{\|h_t-f\| + \|\mu_t-\mu\|_*}$$

$$= \limsup_{t \downarrow 0} \frac{\langle \lambda, tf \rangle - \langle t\mu, f \rangle}{\|h_t-f\| + \|\mu_t-\mu\|_*}$$

$$= \limsup_{t \downarrow 0} \frac{t \sum_{j=1}^m \alpha_j u(s_j) f(s_j) - t\|f\|^2}{\|h_t-f\| + \|\mu_t-\mu\|_*}$$

$$= \frac{\|u\|\|f\| - \|f\|^2}{2\|f\|}$$

$$= \frac{\|u\| - \|f\|}{2} > 0.$$

By (5.21), this implies $\lambda \notin \widehat{D}^* J(f,\mu)(f^{**})$. □

**Corollary 5.7.** *Let $f, u \in C^+[0, 1]$ with $f$ induces $f^{**} \in rca^*[0, 1]$. Suppose that both $f$ and $u$ are increasing with $u(1) > f(1) > 0$. Define $\mu, \lambda \in rca[0, 1]$ by*

$$\mu(1) = f(1) \quad \text{and} \quad \mu(E) = 0, \text{ for any } E \subseteq [0,1);$$

$$\lambda(1) = u(1) \quad \text{and} \quad \lambda(E) = 0, \text{ for any } E \subseteq [0,1).$$

*Then, $\mu \in J(f)$ and $\lambda \in J(u)$ that satisfy*

$$\lambda \notin \widehat{D}^* J(f,\mu)(f^{**}).$$

*Proof.* Since $f, u \in C^+[0, 1]$ and both $f$ and $u$ are increasing, this implies $f(1) = \|f\|$ and $u(1) = \|u\|$. Hence, $1 \in M(u) \cap M(f)$. Rest of the proof is similar to the proof of Theorem 5.6 and it is omitted here. □

**Theorem 5.8.** *Let $f \in C^+[0, 1]$ with $\|f\| > 0$ that induces $f^{**} \in rca^*[0, 1]$. Then, for any $\mu \in J(f)$,*

$$c\mu \notin \widehat{D}^* J(f,\mu)(f^{**}), \text{ for any } c > 0 \text{ with } c \neq 1.$$

*Were, $(c\mu)(A) := c\mu(A)$, for any $A \in \Sigma$.*

*Proof.* Case 1. $c > 1$. For any given $\mu \in J(f)$, we have

$$c\mu \in \widehat{D}^* J(f,\mu)(f^{**}) \quad \Leftrightarrow \quad \limsup_{\substack{(g,\gamma) \to (f,\mu) \\ g \in C[0,1] \text{ and } \gamma \in J(g)}} \frac{\langle c\mu, g-f \rangle - \langle f^{**}, \gamma-\mu \rangle_*}{\|g-f\| + \|\gamma-\mu\|_*} \leq 0$$

$$\Leftrightarrow \quad \limsup_{\substack{(g,\gamma) \to (f,\mu) \\ g \in C[0,1] \text{ and } \gamma \in J(g)}} \frac{c\langle \mu, g-f \rangle - \langle \gamma-\mu, f \rangle}{\|g-f\| + \|\gamma-\mu\|_*} \leq 0. \tag{5.23}$$

For any given $t > 0$, we define $h_t$ by

$$h_t(s) = (1 + t)f(s), \text{ for all } s \in [0,1].$$

Then, $\mu_t := (1 + t)\mu \in J(h_t)$, for any $t > 0$. Similar to the proof of the previous theorems, one can show that

$$(h_t, \mu_t) \to (f, \mu), \text{ as } t \downarrow 0. \tag{5.24}$$

By (5.24), calculating the limit in (5.23), we have

$$\limsup_{\substack{(g,\gamma) \to (f,\mu) \\ g \in C[0,1] \text{ and } \gamma \in J(g)}} \frac{c\langle \mu, g-f \rangle - \langle \gamma - \mu, f \rangle}{\|g - f\| + \|\gamma - \mu\|_*}$$

$$\geq \limsup_{(h_t, \mu_t) \to (f,\mu)} \frac{c\langle \mu, h_t - f \rangle - \langle \mu_t - \mu, f \rangle}{\|h_t - f\| + \|\mu_t - \mu\|_*}$$

$$= \limsup_{t \downarrow 0} \frac{c\langle \mu, tf \rangle - \langle t\mu, f \rangle}{\|h_t - f\| + \|\mu_t - \mu\|_*}$$

$$= \frac{(c-1)\|f\|^2}{2\|f\|}$$

$$= \frac{(c-1)\|f\|}{2} > 0.$$

By (5.23), this implies $c\mu \notin \widehat{D}^* J(f, \mu)(f^{**})$.

Case 2. $0 < c < 1$. In this case, for every $0 < t < 1$, we define

$$h_t(s) = (1 - t)f(s), \text{ for all } s \in [0,1].$$

Similar to the proof of case 1, one has

$$\limsup_{t \downarrow 0} \frac{c\langle \mu, tf \rangle - \langle t\mu, f \rangle}{\|h_t - f\| + \|\mu_t - \mu\|_*}$$

$$\geq \frac{(1-c)\|f\|^2}{2\|f\|}$$

$$= \frac{(1-c)\|f\|}{2} > 0. \qquad \square$$

**Acknowledgments.** The author is very grateful to Professor Boris S. Mordukhovich, Professor Simeon Reich and Professor Christiane Tammer for their kind communications, valuable suggestions and enthusiasm encouragements in the development stage of this paper.

**Appendix**

**Some properties of the normalized duality mapping in Banach spaces**

The inverse of the normalized duality mapping $J$ is denoted by $J^{-1}$, which is a set-valued mapping from $X^*$ to $X$ such that, for any $\varphi \in X^*$, we have

$$J^{-1}(\varphi) = \{x \in X : \varphi \in J(x)\}.$$

Then, for any $x \in J^{-1}(\varphi)$, we have $\|x\| = \|\varphi\|_*$. For $x, y \in X$, if $J(x) \cap J(x) \neq \emptyset$, then $x$ and $y$ are said to be generalized identical. For $x \in X$, the set of all its generalized identical points is denoted by $\mathfrak{J}(x)$. That is,

$$\mathfrak{J}(x) = \{y \in X : J(x) \cap J(x) \neq \emptyset\}.$$

The generalized identity has the following properties.

(a) $x \in \mathfrak{J}(x)$, for any $x \in X$;
(b) $\|y\| = \|x\|$, for any $y \in \mathfrak{J}(x)$.

We list some properties of the normalized duality mapping below for easy reference. For more details, one may see Sections 4.2–4.3 and Problem set 4.2 in [28] and [12, 17, 23, 25].

($J_1$). For any $x \in X$, $J(x)$ is nonempty, bounded, closed and convex subset of $X^*$;
($J_2$). $J$ is the identity operator in any Hilbert space $H$;
($J_3$). $J(\theta) = \theta^*$ and $J^{-1}(\theta^*) = \theta$;
($J_4$). For any $x \in X$ and any real number $\alpha$, $J(\alpha x) = \alpha J(x)$;
($J_5$). For any $x, y \in X$, $j(x) \in J(x)$ and $j(y) \in J(y)$, $\langle j(x) - j(y), x - y \rangle \geq 0$;
($J_6$). For any $x, y \in X$, $j(x) \in J(x)$, and $j(y) \in J(y)$, we have

$$2\langle j(y), x - y \rangle \leq \|x\|^2 - \|y\|^2 \leq 2\langle j(x), x - y \rangle;$$

($J_7$). $X$ is strictly convex if and only if $J$ is one-to-one, that is, for any $x, y \in X$,

$$x \neq y \implies J(x) \cap J(y) = \emptyset;$$

($J_8$). $X$ is strictly convex if and only if, for any $x, y \in X$ with $x \neq y$,

$$j(x) \in J(x) \text{ and } j(y) \in J(y) \implies \langle j(x) - j(y), x - y \rangle > 0;$$

($J_9$). $X$ is strictly convex if and only if, for any $x, y \in X$ with $\|x\| = \|y\| = 1$ and $x \neq y$,

$$\varphi \in J(x) \implies 1 - \langle \varphi, y \rangle > 0;$$

($J_{10}$). If $X^*$ is strictly convex, then $J$ is a single-valued mapping;
($J_{11}$). If $X$ is reflexive, then $X$ is smooth if and only if $J$ is a single-valued mapping;
($J_{12}$). If $J$ is a single-valued mapping, then $J$ is a norm to weak$^*$ continuous;
($J_{13}$). $X$ is reflexive if and only if $J$ is a mapping of $X$ onto $X^*$;

($J_{14}$). If $X$ is smooth, then $J$ is a continuous operator;
($J_{15}$). $J$ is uniformly continuous on each bounded set in uniformly smooth Banach spaces;
($J_{16}$) $X$ is strictly convex, if and only if $\Im(x) = x$, for every $x \in X$.